\documentclass[]{article}

\usepackage[utf8]{inputenc}
\usepackage{graphicx} % Required for inserting images
\usepackage{amsthm}
\usepackage{amssymb}
\usepackage[UKenglish]{babel}
\usepackage{colonequals}
\usepackage{mathtools}
\usepackage{xcolor}
\usepackage{upgreek}

\usepackage{amsmath}	% Mathe: Für align-Umgebung
\usepackage{amsfonts}   % Mathe für mathbb
\usepackage{enumitem}   % Custom enumerate

\usepackage{amsthm}
\usepackage[hidelinks]{hyperref}
\usepackage{tikz}
\usepackage{tikz-cd}    % Kommutierende Diagramme
\usetikzlibrary{arrows, calc, angles, quotes, 3d}
\usepackage{comment}
\usepackage{mathrsfs}
\usepackage{makeidx}

\usepackage{tabularray}
\UseTblrLibrary{amsmath}

\title{Promise Graph Coloring is NP-Hard}
\author{Sebastian Meyer}

\newtheorem{theorem}{Theorem}[section]
\newtheorem{proposition}[theorem]{Proposition}
\newtheorem{lemma}[theorem]{Lemma}

\newtheorem{corollary}[theorem]{Corollary}

\newtheorem{conjecture}[theorem]{Conjecture}

\theoremstyle{definition}
	\newtheorem{definition}[theorem]{Definition}

	\newtheorem{algdef}[theorem]{Algorithmic Definition}
	\newtheorem{problem}{Problem}[section] %Used for Decision problems
	\newtheorem{claim}[theorem]{Claim}
\theoremstyle{remark}
	\newtheorem{remark}[theorem]{Remark}

\newenvironment{proofof}[1]{\begin{proof}[Proof of #1]}{\end{proof}}

\renewcommand{\subset}{\subseteq}

\newcommand{\mitte}{\,\middle|\,}

\newcommand{\de}[1]{\emph{#1}\index{#1}}

\newcommand{\id}{\operatorname{id}}
\newcommand{\iso}{\cong}	%Isomorphismus

\newcommand{\Hom}{\operatorname{Hom}}	% Homomorphisms in a category
\newcommand{\map}{\operatorname{map}}	% set of all maps

\newcommand{\blank}{-}	%Freie Stelle in Multimorphismen
\newcommand{\IN}{\mathbb{N}}
\newcommand{\IZ}{\mathbb{Z}}
\newcommand{\IR}{\mathbb{R}}

\newcommand{\IQ}{\mathbb{Q}}

\newcommand{\powerset}{\mathcal{P}}

\newcommand{\ts}[1]{t}%\operatorname{ts}_{#1}}

\newcommand{\ix} {\fdual{\mathscr{X}}}
\newcommand{\iix}{\cdual{\mathscr{X}}}

\newcommand{\inst}{\mathcal{I}}

\newcommand{\deld}{\,\textup{d}}	% Integrations-d
\newcommand{\dist}{\operatorname{dist}}
\newcommand{\distloc}{\operatorname{dist}_{\text{loc}}}
\newcommand{\cdual}{\widehat}	% Fourier Dual
\newcommand{\fdual}{\tilde}	% Finite Fourier Dual
\newcommand{\IT}{\mathbb{T}}	% Torus

\newcommand{\colim}{\operatorname{colim}}

\newcommand{\IE}{\mathbb{E}}
\newcommand{\lebesque}{\mathcal{L}}

\begin{document}

\maketitle

\begin{abstract}
	Given any pair $(G,H)$ of finite graphs (undirected, without loops, non-bipartite) with a graph homomorphism $G\to H$, the decision problem that decides for any graph $K$ wether there is a graph homomorphism $K\to G$ or no graph homomorphism $K\to H$ (where it is promised that one case happens) is NP hard.
	In other words, the Brakensiek-Guruswami-Conjecture is true.
	
	In addition to the new result, this paper should serve as a more readable version of OpenAIs paper 106 on hardness of finding large independent sets
	in three-colorable graphs.
\end{abstract}
\tableofcontents
\subsection*{Acknowledgments}
I thank Libor Barto and Jakub Opr\v{s}al for helpful remarks on theorems I wanted to use and I thank Antoine Motettet for pointing out the existing parallel work.

\subsection*{AI disclosure}
No AI was used to obtained the result. It is a generalization of \cite{OAI_PCSP} by OpenAI which is properly cited.

\subsection*{Parallel Work}
A proof to theorem was completely independently published in the preprints \cite{hexagon:2611.00619v1}, ten days ago and in \cite{hexagon:2611.01095v1} two days ago.

\section{Introduction}
The typical setting of graph coloring is given any graph $K$, decide wether the vertices of $K$ can be colored in $k$ colors such that no neighboring vertices admit the same color. This problem is known to be NP hard if and only if $k\ge 3$. We generalize this question in two directions:
\begin{enumerate}
	\item Instead of a coloring in equal colors, we color the vertices of $K$ in the vertices of another graph $G$ where for two vertices $(v,w)$ in $K$ that are connected by an edge, also their images in $G$ have to be connected by an edge. This problem was firstly classified in \cite{HN90}. It can be solved in polynomial time if and only if $G$ is bipartite or has a loop.
	\item Instead of coloring in a single number $k$ of colors, we introduce a promise: For two numbers $k_1\le k_2$, we ask wether any graph $K$ can be colored in $k_2$ or not even in $k_1$ colors where it is promised that one of the two cases occurs. This problem was recently solved in \cite[Corollary 1.7]{Hardness4to1games}.
\end{enumerate}
We consider a unified generalization. A graph homomorphism $G$ to $H$ is a map $f$ from the vertices of $G$ to the vertices of $H$ such that $(f(v),f(w))$ is an edge whenever $(v,w)$ is one.
\begin{problem}[Promise Graph Coloring]
	Parameter are two (finite, undirected) graphs $(G,H)$ such that there is a graph homomorphism $G\to H$.
	\begin{description}
		\item[Input:] A graph $K$
		\item[Output: Yes] if there is a graph homomorphism $K\to G$
		\item[Output: No] if there is no graph homomorphism $K\to H$
	\end{description}
	In other cases, there can be any output.
\end{problem}

We show that even in this case, the problem is always NP hard, except for trivial cases:
\begin{theorem}[Brakensiek-Guruswami-Conjecture, Conjectured in \cite{BG18}] \label{TheoremMain}
	The promise graph coloring problem is NP complete whenever $G$ and $H$ are non-bipartite and have no loop and in P otherwise.
\end{theorem}

\subsection{Proof Overview}\label{Section_overview}
The proof generalizes the recent paper \cite{OAI_PCSP} by OpenAI. We give a full proof of the generalization and capture all results from \cite{OAI_PCSP} as the original paper is very hard to read.

We introduce gap label cover and promise $(C_k,K_k)$ graph coloring in Section~\ref{Section_Preliminaries} and give the easy reductions. Additionally, we state results we use from other papers, namely NP-hardness of the gap label cover problem, a Ramsey theory result and an analytic version of the Junta theorem.

In Section~\ref{Section_Reduction}, we give the main reduction step by describing the reduction from promise $(C_k,K_k)$ graph coloring to $\sigma$-gap $A$-multi-label cover. We also include a description of how to compute the constants $\sigma$ and $A$ but we will not compute them.

The correctness of the algorithm is shown in Section~\ref{Section_Correctness}.

\begin{figure}
	$$
	\begin{tikzpicture}
		\node [draw, align=center, text width=120] (A) at (0,3) {$\forall \, G,H\text{ non-bipartite}\colon $\\ Promise $(G,H)$ Graph Coloring\\ is NP-hard};
		\node [draw, align=center, text width=120] (B) at (0,0) {$\forall \, k\ge 3 \text{ odd}\colon$\\Promise $(C_k,K_k)$ Graph Coloring\\ is NP-hard};	
		\node [draw, align=center, text width=120] (C) at (-7,3) {$\exists\, \sigma, B\colon$\\$\sigma$-Gap Label Cover with bound $B$\\ is NP-hard};
		\node [draw, align=center, text width=120] (D) at (-7,0) {$\forall \, \sigma, A\  \exists \, B\colon$\\$\sigma$-Gap $A$-Multi-Label Cover with bound $B$\\ is NP-hard};	
		\draw[->] (A) edge (B) (B) edge node[above]{hard part} (D) (D) edge (C);
	\end{tikzpicture}
	$$
	
	\caption{General Strategy}
	\label{Figure_Overview}
\end{figure}
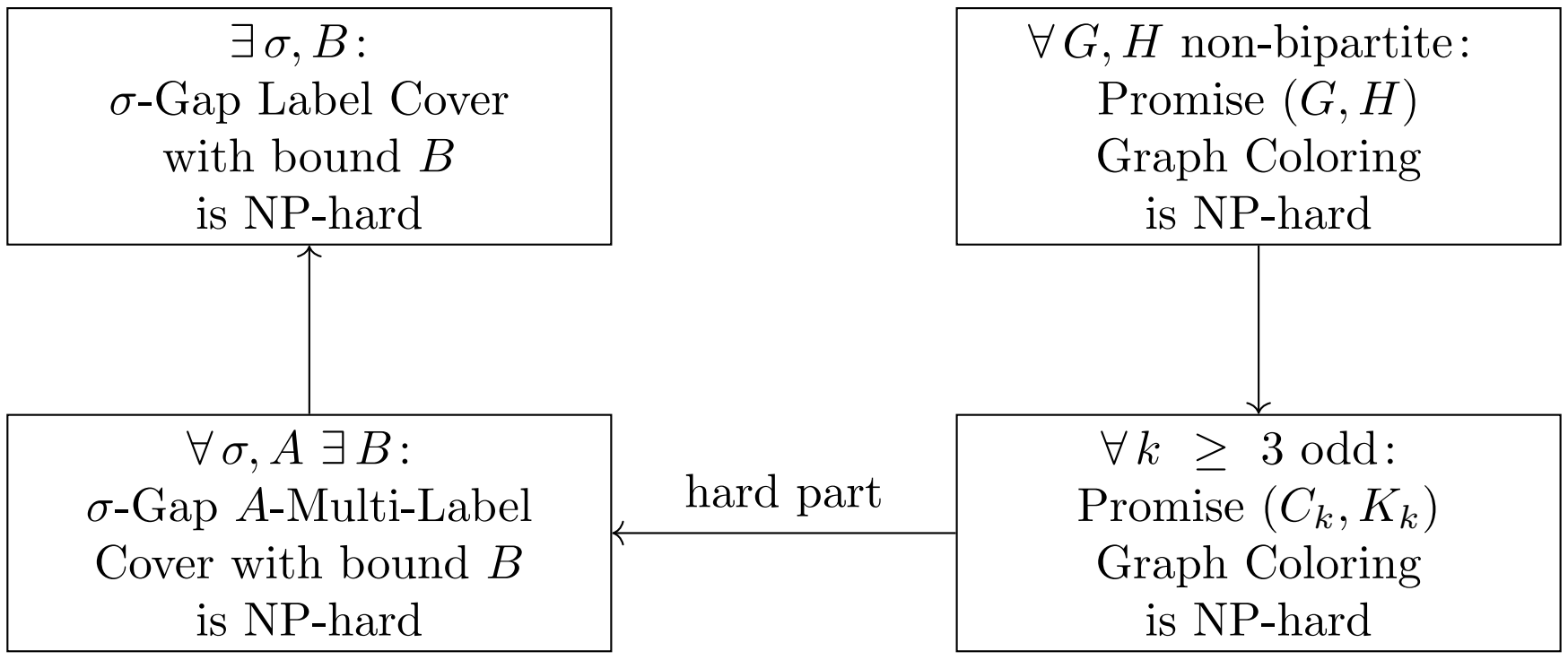

Mathematically, we reduce the promise graph coloring problem on the gap label cover problem in three steps as shown in Figure~\ref{Figure_Overview}: At first we recall the well-known fact that we may assume the two graphs $(G,H)$ to be $C_k$, an odd cycle on $k$ vertices and $K_k$, the full graph on $k$ vertices. In a second step, we generalize the gap label cover problem to gap multi-label cover where multiple labels are allowed instead of just one.

Finally, we reduce the promise $(C_k,K_k)$-problem (where $k$ is odd) to the gap multi-label cover problem:
\begin{itemize}
	\item We associate to every instance of the gap multi-label cover problem a graph $G$, which is essentially given by the Pontryagin dual or finite Fourier transform of the label cover instance (Section~\ref{Section_Reduction}).
	\item If the label cover instance has a solution then evaluating induces a map from the Pontryagin dual to the unit sphere, which easily translates to a map from $G$ to $C_k$ (Section~\ref{Section_corr_SatInst}).
	\item To show the other direction, we may assume to the contrary that the gap multi-label cover instance has no approximate solution and that there is a homomorphism $G \to K_k$.
	\item We define for all constraint tuples $\pi_c$ in the label cover instance (strictly speaking in a more complicated input $w\in W_-$) a probability measure $\nu_w$ on a modified version of $G$ that is almost $G^r$. By averaging over all tuples, we get a total probability measure $\nu$ on $G^r$ (Section~\ref{Section_corr_ProbDist}).
	\item If $G\to K_k$, then $G^r\to K_k$ and there is an independent vertex set $\mathcal A\subset G^r$ with measure $\nu(\mathcal A)\ge 1/k$ (start of Section~\ref{Section_corr_IndependentSet}).
	\item As $G$ is a Pontryagin dual, $G$ and also $G^r$ carry a metric. We define a Lipschitz continuous map $f\colon G^r \to [-1,1]$ that is 0 on average and 1 on $\mathcal A$ (end of Section~\ref{Section_corr_IndependentSet}).
	\item Applying a version of the Junta theorem to $f$ implies that $f$ essentially equals a map $g$ that only depends on a finite number $d$ of inputs. These inputs define a multi-label cover for the original instance (Section~\ref{Section_corr_JuntaApprox}).
	\item By assumption, this multi-label cover does not has an approximate solution, so many of these coordinates are not preserved by the label cover maps $\pi_c$ (Section~\ref{Section_corr_NegInst}).
	\item Technically, we want to use the previous statement for $w$ and $G^r$ rather than for $\pi_c$ and $G$ which introduces another technicality which we can bound by a Ramsey argument such that it occurs only for few $w$ (Section~\ref{Section_corr_Ramsey}).
	\item For every other $w$ (which roughly correlates to a map $\pi_c$ that does not preserves the input coordinates of $g$) the map $g$ is close to zero, measured against $\nu_w$. Thus, by some analytic approximations, also $f$ is close to zero (Section~\ref{Section_corr_Analysis}).
	\item By combining the previous estimates and choosing quantifiers in the right order, we get that 
	$$
	\nu(\mathcal A)
	= \sum_{w\in W_-}\frac{1}{|W_-|} \nu_w\{\mathcal A\}
	\le \sum_{w\in W_-}\frac{1}{|W_-|} \sum_{v\in G} f(v)^2 \nu_w(\{v\})
	$$
	becomes arbitrarily small as for many $w$ the value of $f(v)^2$ measured against $\nu_w$ is small. We may even assume that this value is smaller than $1/k$, a contradiction (Section~\ref{Section_corr_Combined}).
\end{itemize}

\subsection{Some Previous Results}
Many special cases of the promise graph coloring problem have been proven before.
The graph coloring problem without any promise was classified in \cite{HN90} and again in \cite{Bul05}. The promise graph coloring 3 coloring 4 colorable graphs was shown to be NP-hard in \cite{KLS00} and \cite{GK04}. Hardness of 3 coloring 5 colorable graphs was proven in \cite{BBKO21}. 

The conjecture that classifies promise graph coloring for all finite undirected graphs was stated in \cite{BG18}.
Below the number 3, there was also some progress: Hardness of $(2+\varepsilon)$-coloring $(4-\varepsilon)$-colorable graphs for any $\varepsilon>0$ was shown in \cite{KOWZ23} and even hardness of $(2+\varepsilon)$-coloring 4-colorable graphs in \cite{AFOTW25}.
Hardness of 3-coloring $n$ colorable graphs for any $n\ge 3$ was proven in \cite{Hardness4to1games} and shortly thereafter in \cite{OAI_PCSP}. We generalize the proof from \cite{OAI_PCSP} to hardness of $(2+\varepsilon)$-coloring $n$-colorable graphs and repeat the well-known proof that this suffices to completely classify promise graph coloring. 

After finishing this draft, the author was pointed towards the fact that during the writing, new progress has happened. 
This theorem was completely independently shown in the preprint \cite{hexagon:2611.00619v1} by Yuwen Zhou through an interactive dialogue with GPT-6.2 Astra (OpenAI) and Claude Opus 5.5 (Anthropic), less than two weeks ago. The preprint \cite{hexagon:2611.00619v1} seems to use a similar approach and sometimes even the same variable names as \cite{OAI_PCSP} and this paper.

Moreover, two days ago, another independent proof of the statement was put out in \cite[Theorem 2.1]{hexagon:2611.01095v1} by Brakensiek, Guruswami, Kozik and Ren with the use of ChatGPT-6 Astra. It uses slightly different ideas and avoids the PCP theorem.

\section{Preliminaries}
\label{Section_Preliminaries}
\subsection{General Notions}
We denote the set $\{1,2,\dots, n\}=[1,n]\cap \IZ$ by $[n]$.

For an integer $D$, we denote the group $D^{-1}\IZ/\IZ$ by $\IT_D$ and consider it as part of the full circle $\IR/\IZ$, which we denote by $\IT$ and consider to be both a group with addition and a metric space where $\dist(a+\IZ,b+\IZ)$ is defined as minimum $\min_{n\in \IZ}(a-b+n)$. (The group $D^{-1}\IZ/\IZ$ is isomorphic to $\IZ/(D\IZ)$, but the metric is different.)
We denote the addition by $1/2$ as map $\IT \to \IT$ respectively $\IT_D\to \IT_D$ when $D$ is even with $T$. 
We denote the set of maps from a set $X$ to a set $Y$ by $\map(X,Y)$.

\subsection{Graph Coloring}
For a graph $G$, we denote its vertex set with $V(G)$ or simply $G$ and its edge set with $E(G)$. We assume all graphs to be finite and have no parallel edges. A \de{graph homomorphism} from $G$ to $H$ is a map $f$ from $V(G)$ to $V(H)$ such that $\{v_1,v_2\}\in E(G)\implies \{f(v_1),f(v_2)\}\in E(H)$. We write $G\to H$ if such a homomorphism exists. The product $G\times H$ of the graphs $G$ and $H$ has vertex set $V(G)\times V(H)$ and an edge from $(v_1,v_2)$ to $(v_3,v_4)$ if and only if $(v_1,v_3)$ is an edge in $G$ and $(v_2,v_4)$ an edge in $H$.
By $K_k$, we denote the complete graph on the vertex set $[k]\coloneqq \{1,2,3,\dots,k\}$. By $C_k$, we denote the odd cycle on $k$ vertices, that is the graph
$$
(V(C_k),E(C_k))=([k],\{\{1,2\},\{2,3\},\dots,\{k-1,k\},\{k,1\}\}).
$$
More general, we define for positive integers $a,b$ the graph $K_\frac{a}{b}$ as
$$
(V(K_\frac{a}{b}),E(K_\frac{a}{b}))=(a^{-1}\IZ/\IZ,\{\{v_1,v_2\} \mid \dist(v_1-v_2)\ge b/a \})
$$
and note that it has no edges for $a/b<2$.
\begin{remark} The following holds. \label{Remark_OrderCompleteGraphs}
	\begin{align*}
		\forall k\in \IN:{}&& K_\frac{k}{1} &\iso K_k\\
		\forall k\in \IN:{}&& K_\frac{2k+1}{k} &\iso C_{2k+1}\\
		\forall a,b,c,d\in \IN:{}&a/b\le c/d \implies &K_\frac{a}{b} &\to K_\frac{c}{d}
	\end{align*}
\end{remark}

\begin{lemma}
	Consider two finite graphs $G$ and $H$ with $G\to H$, $G$ non-bipartite and $H$ without loop.
	Then, the promise $(G,H)$-graph coloring problem reduces in constant (and thus polynomial) time to the promise $(C_k,K_k)$-graph coloring problem where $k$ is any odd integer larger than the size of $G$ and $H$ by taking the same inputs.
\end{lemma}
\begin{proof}
	As $G$ is non-bipartite, it embeds an odd cycle $C_{2n+1} \to G$ with $2n+1\le |V(G)|$. As $H$ has no loop, choosing any bijection $V(H) \to [|V(H)|]$ induces a graph homomorphism $H\to K_{|V(H)|}$. Now, Remark~\ref{Remark_OrderCompleteGraphs} implies
	$$
	C_k \iso K_\frac{k}{2^{-1}(k-1)} \to K_\frac{2n+1}{n} \iso C_{2n+1} \to G \to H \to K_{|V(H)|} \to K_k
	$$
	and thus any input graph that maps to $H$ also maps to $K_k$ and any input graph that does not map to $G$ also does not map to $C_k$.
\end{proof}

\subsection{Gap Label Cover}

Gap label cover is a promise decision problem.

\begin{problem}[$\sigma$-Gap Label Cover Problem with Bound $B$]
	Depends on parameters $\sigma \in (0;1]$ and $B\in \IN, B\ge 2$.
	\begin{description}
		\item[Input:] A label cover instance with bound $B$, that are 
		\begin{itemize}
			\item three numbers $n,m,C\in \IN$,
			\item a list $(U_1,\dots,U_n)$ of sets of size at most $B$,
			\item a list $(U_{n+1},\dots,U_{m+n})$ of sets of size at most $B$,
			\item For all $c\in [C]$ numbers $i_c\in [n]$, $j_c\in [m+n]\setminus[n]$ and a map $\pi_c\colon U_{i_c} \to V_{j_c}$.
		\end{itemize}
		
		\item[Output: Yes] if there are points $u_i\in U_i$ for all $i\le m+n$ with $\pi_c(u_{i_c})=v_{j_c}$.
		\item[Output: No] if for every points $(u_i)_{i\in [m+n]}\in U_i^{n+m}$ the number of $c\in [C]$ with $\pi_c(u_{i_c})=v_{j_c}$ is less than $\delta C$.
	\end{description}
	In other cases, there can be any output.
\end{problem}

We also consider slight modifications. An infinite label cover instance allows the sets $U_1,\dots,U_{n+m}$ to be infinite.
\begin{problem}[$\sigma$-Gap $A$-Multi-Label Cover Problem with Bound $B$]
	Depends on parameters $\sigma \in (0;1]$, $A\in \IN \setminus\{0\}$ and $B\in \IN, B\ge 2$.
	\begin{description}
		\item[Input:] A label cover instance with bound $B$.
		\item[Output: Yes] if there are points $u_i\in U_i$ for all $i\le m+n$ with $\pi_c(u_{i_c})=v_{j_c}$.
		\item[Output: No] if for every tuple of subsets $U'_i\subseteq U_i$ (where $i$ is from $[m+n]$), of size $|U'_i|\le A$, the number of $c\in [C]$ such that there is $u\in U'_{i_c}$ with $\pi_c(u)\in U'_{j_c}$ is less than $\delta C$.
	\end{description}
	In other cases, there can be any output.
\end{problem}

The PCP theorem \cite{PCPThm} implies the following result for gap label cover.
\begin{theorem}\label{Theorem_PCP}
	There are $\sigma<1$ and $B \in \IN$ such that the $\sigma$-gap label cover problem with bound $B$ is NP-hard.
\end{theorem}
Already $B=7$ works by translating approximate 3-SAT to label cover instances. The statement for approximate 3-SAT is shown in \cite[Theorem 1.2.2]{PCPThm}.

\begin{lemma} \label{Lemma_LabelreductionInduction}
	The $\sigma$-gap label cover problem with bound $B$ reduces in polynomial time to the $\sigma^2$-gap label cover problem with bound $B^2$.
\end{lemma}
\begin{proof}
	Given any input for the $\sigma$-gap label cover problem with bound $B$, consider the following input for the $\sigma^2$-gap label cover problem with bound $B^2$:
	\begin{itemize}
		\item $n^2, m^2, C^2$
		\item $(U_{i} \times U_{j})_{i,j\in [n]}$ (where we identify $[n] \times [n]$ with $[n^2]$)
		\item $(U_{i} \times U_{j})_{i,j\in [m+n]\setminus [n]}$ (with a similar identification)
		\item For all $(c,e)\in [C]^2 \iso [C^2]$ define $i_{c,e}=(i_c,i_e)$, $j_{c,e}=(j_c,j_e)$ and $\pi_{c,e}(u,w)=(\pi_c(u),\pi_c(w))$.
	\end{itemize}
	It is correct and computable in polynomial time.
\end{proof}

\begin{lemma} \label{Lemma_Multilabelreduction}
	The $\sigma$-gap $A$-multi-label cover problem with bound $B$ reduces in constant time to the $\sigma\cdot A^2$-gap label cover problem with bound $B$ by taking the same input.
\end{lemma}
\begin{proof}
	Choose $u_i$ at random inside $U'_i$ for all $i$. Then, the expected value for each satisfiable map $\pi_c$ to hit a satisfiable tuple is at least $1/A^2$. Thus the average number of satisfied clauses is $\sigma\cdot A^2$ and therefore there is a case with at least $\sigma\cdot A^2$ satisfied clauses.
\end{proof}

\begin{corollary} \label{Corollary_BoundLabelCover}
	For every $\sigma\in (0,1), A\in \IN \setminus\{0\}$, there exists $B \in \IN$ such that the $\sigma$-gap $A$-multi-label cover problem with bound $B$ is NP-hard.
\end{corollary}
\begin{proof}
	Apply Lemma~\ref{Lemma_Multilabelreduction} and then inductively Lemma~\ref{Lemma_LabelreductionInduction} until the current $\sigma$ is smaller than the $\sigma$ from Theorem~\ref{Theorem_PCP}. As enlarging $\sigma$ just makes the claim stronger, Theorem~\ref{Theorem_PCP} can be applied to obtain hardness.
\end{proof}
\begin{remark} \label{Remark_BoundLabelCoverQuantitative}
	By quantifying the proof, the constant $B$ in Corollary~\ref{Corollary_BoundLabelCover} can be chosen to be the smallest integer larger than $$\exp\left(2\frac{\log(B_0)\log(\sigma d^2)}{\log(\sigma_0)}\right)$$
	where $B_0$ and $\sigma_0$ are the constants from Theorem~\ref{Theorem_PCP}.
\end{remark}
\begin{remark}
	We may assume that the maps $\pi_c\colon U_{i_c} \to V_{j_c}$ of a label cover instance are surjective as we can reduce the general case to the surjective case by repeatedly removing non-surjective points, which is a polynomial time algorithm.
\end{remark}
\subsection{Junta Theorem for Graphs}
To show that the reduction from gap label cover to graph coloring, which we will define in Section~\ref{Section_Reduction}, is correct, we will need a variation of the Junta Theorem for graphs. We state it here for reference but we will not use it.
\begin{theorem}[Theorem 1.1 of \cite{GraphJuntas}]
	Let $G$ be a connected non-bipartite graph without loops and denote $V=V(G), E=E(G)$. Then, there exists a function $j\colon \IR \to \IN$ such that if $I\subseteq V^n$ is an independent set in $G^n$ then for every $\varepsilon>0$ there exists a set $J\subseteq V^n$ depending on at most  $j(\varepsilon)$ coordinates such that $\mu(I\setminus J)\le \varepsilon$ and such that $J$ spans less than $\varepsilon|E(G^j)|$ edges in the graph $G^j$.
	
	Here, $\mu$ is the product measure on $V^n$ of the probability measure on $V$ that associates to a vertex the number of incoming edges divided by $2|E|$.
\end{theorem}

\subsection{Analytic Junta Theorem}
We actually use an analytic version of the Junta Theorem and some basic analytic definitions that we introduce here.

\begin{definition}
	A \de{pseudo-metric} space $X$ is a set with a map $\dist\colon X^2\to [0,\infty]$ such that 
	\begin{align*}
		\forall x\in X\colon && \dist(x,x)&=0 \\
		\forall x,y\in X\colon && \dist(x,y)&=\dist(y,x)\\
		\forall x,y,z\in X\colon && \dist(x,y)+\dist(y,z)&\ge \dist(x,z)
	\end{align*}
	holds.
	
	A map $f$ from a pseudo-metric space $X$ to $\IR\cup \{\infty\}$ is called \de{$L$-Lipschitz continuous}, if 
	$$
	\forall x,y \in X: \dist(x,y)<\infty \implies f(x)-f(y)\le L \dist(x,y)
	$$
	holds.
\end{definition}
Note that every Lipschitz continuous map is continuous and uniformly continuous.

\begin{theorem}[Analytic Junta Theorem for $L^1$, c.f. Theorem 1.1 of \cite{AnalyticJunta}]
	\label{TheoremJunta2a}
	
	For every $0\le L<\infty$ and $\tau>0$ there is an integer $J(L,\tau)\ge 1$ such that for every positive integer $N$, every $L$-Lipschitz function $f\colon [0,1]^N \to [-1,1]$ can be approximated by a function $g$ that depends on at most $J(L,\tau)$ coordinates such that
	$$
	\int_{[0,1]^N} |f(x)-g(x)| \deld x <\varepsilon
	$$
	holds.
\end{theorem}
Note that $J(L,\tau)$ does not depend on $N$.

The theorem in \cite[Theorem 1.1]{AnalyticJunta} discusses the more general case of any connected locally connected metric space and uniformly continuous functions instead of just Lipschitz continuous ones.

We can also show a quadratic bound.
\begin{theorem}[Analytic Junta Theorem for $L^2$, Theorem 2.2 of \cite{OAI_PCSP}]
	\label{TheoremJunta2}
	
	For every $0\le L<\infty$ and $\tau>0$ there is an integer $J(L,\tau)\ge 1$ such that for every positive integer $N$, every $L$-Lipschitz function $f\colon [0,1]^N \to [-1,1]$ can be approximated by a function $g$ that depends on at most $J(L,\tau)$ coordinates such that
	$$
		\int_{[0,1]^N} (f(x)-g(x))^2 \deld x <\tau^2
	$$
	holds.
\end{theorem}
\begin{proof}
	Theorem~\ref{TheoremJunta2a} just discusses the $L^1$-approximation, so applying it gives for every $\varepsilon$ a function $g$ with
	$$
	\int_{[0,1]^N} |f(x)-g(x)| \deld x <\varepsilon
	$$
	As we may assume that $f(x)$ and $g(x)$ map into $[-1,1]$, we get $|f(x)-g(x)|<2$ and thus 
	$$
	\int_{[0,1]^N} (f(x)-g(x))^2 \deld x <2\varepsilon\le \tau^2
	$$
	holds when choosing $\varepsilon\le \tau^2/2$.
\end{proof}

As taking the integral respectively average over some or all inputs will occur multiple times, we give it its own notion. For any function $f\colon X_1\times \dots \times X_n \to Y$ and any $i\in [n]$, denote by $\IE_i f$ the map 
\begin{align*}
	X_1\times \dots \times X_n &\to Y \\
	(x_1,\dots,x_n) &\mapsto \int_{X_i} f(x_1,\dots,x_{i-1}, y,x_{i+1},\dots,x_n ) \deld \mu(y)
\end{align*}
which implicitly depends on a probability measure $\mu$ on $X_i$. The probability measure on $[0,1]$ and on $\IT=\IR/\IZ$ is the Lebesque measure unless stated otherwise. 

In addition to the main application of the Theorem, we will use it a second time together with a Ramsey argument to approximate an error term. The second application uses the theorem in a slightly different version:
\begin{theorem}[Lemma 2.3 of \cite{OAI_PCSP}] \label{TheoremJunta3}
	For every $0\le L<\infty$ and $\varepsilon>0$, there are $\beta>0$ and an integer $K\ge 1$ such that for all $s\in \IN$, every $L$-Lipschitz function $H\colon \IT^s \to [-1,1]$ with $\int_{\IT^s} H(z) \deld z = 0$ has this properties:
	\begin{enumerate}
		\item If $\int_{\IT^s} (H(z))^2 \deld z > \varepsilon$ then there is $i\in [s]$ such that $$\int_{\IT^s} (H(x)-\IE_i H(x))^2\deld x\ge \beta^2$$ holds and
		\item at most $K$ indices $i$ satisfy  $\int_{\IT^s} (H(x)-\IE_i H(x))^2\deld x\ge (\beta/2)^2$.
	\end{enumerate}
\end{theorem}
\begin{proof}
	Consider $f\colon [0,1]^s \to [-1,1]$ by $f(x)=H(x\pmod{1})$ and choose an approximation
	$g\colon [0,1]^s \to [-1,1]$ of $f$ according to Theorem~\ref{TheoremJunta2} with bound 
	$$\int_{[0,1]^s} (f(z)-g(z))^2\deld z<\varepsilon/4$$ 
	that depends on a set $T\subset[s]$ of at most $J(L,\sqrt{\varepsilon}/2)$ coordinates. We choose $g$ as a partial arithmetic mean of $f$, so we get $\int_{\IT^s} g(z) \deld z=0$.
	Now, we get assuming $\int_{\IT^s} (H(z))^2 \deld z>\varepsilon$ that
	$$
		\int_{[0,1]^s} g(z)^2\deld z>3\varepsilon/4>\varepsilon/4
	$$
	by triangle inequality. However, an analytic approximation implies
	$$
		\int_{[0,1]^s} g(z)^2\deld z \le \sum_{i\in T}  \int_{\IT^s} (H(z)-\IE_i H(z))^2\deld z
	$$
	so there is $i\in T$ with 
	$$\int_{\IT^s} (H(z)-\IE_i H(z))^2\deld z \ge \varepsilon/(4|T|)=\varepsilon/(4 J(L,\sqrt{\varepsilon}/2))\eqcolon \beta^2$$
	defining $\beta$ this way.
	Now, the first part is constructed to be true and the second part follows easily again by Theorem~\ref{TheoremJunta2} with distance $\tau=\beta$ and $K=J(L,\beta/2)$.
\end{proof}
The proof implies also concrete boundaries:
\begin{remark} \label{RemarkJuntaQuantitative}
	In Theorem~\ref{TheoremJunta3}, $\beta$ can be chosen as $\sqrt{\varepsilon/(4J(L,\sqrt{\varepsilon}/2))}$
	and $K$ can be chosen as $J(L, \beta/2)$ where $J$ is the function from Theorem~\ref{TheoremJunta2}.
\end{remark}

\subsection{Ramsey Theory}
We also use a Lemma that follows by a Ramsey theoretic argument.

The $s$-element subsets of a set $X$ are denoted $\binom{[r]}{s}$, the non-empty subsets by $\powerset^+(X)$ and the non-empty at most $s$-element subsets with $\powerset^+_{\le s}(X)$.

\begin{lemma}[c.f.~Lemma 5.1 of \cite{OAI_PCSP}]
	\label{Lemma_Ramsey}
	For every pair $s,d$ of integers and every real number $\varepsilon$, there exists an integer $r\ge s$ and a probability distribution $\mu$ on $\binom{[r]}{s}$, the $s$-element subsets of $[r]$, with the following property:
	
	For any (possibly infinite) set $M$ and any map $A\colon \powerset^+_{\le s}([r]) \to \powerset^+_{\le d}(M), I\mapsto A_I$ that additionally satisfies
	\begin{align*}
		\forall I,J \in \powerset^+_{\le s}([r])&\colon \max(I)<\min(J) \implies &A_I\cap A_J&=\emptyset.
	\end{align*}
	Then,
	the $\mu$-measure of the set
	$$
		\left\{B\in \binom{[r]}{s} \mitte \forall I,J \in \powerset^+(B)\colon I\cap J=\emptyset \implies A_I\cap A_J=\emptyset\right\}
	$$
	is at least $1-\varepsilon$.
\end{lemma}

Note that $\varepsilon$ does not depend on the map $A$.
We do not provide a proof and refer in the lack of a better source to the appendix of \cite{OAI_PCSP}.

\begin{remark}
	As we are not providing a proof, let us at least give some intuition about the Lemma. Its following variation can easily seen to be equivalent, we are neither proving the variation nor the equivalence.
	
	\begin{theorem}[c.f.~Lemma 5.1 of \cite{OAI_PCSP}]
		For every pair $s,d$ of positive integers $(s\ge 3)$ and every real number $\varepsilon>0$, there exists an integer $r\ge s$ and a probability distribution $\mu$ on $\binom{[r]}{s}$ with the following property:
		
		For every map $A\colon \powerset^+_{\le s}([r]) \to \powerset^+_{\le d}([r])$ with 
		$$
		\forall I \in \binom{[r]}{s-1}\colon \min I \le \min A(I) \le \max A(I) \le \max I,
		$$
		the inequality
		\begin{equation*}
			\mu\left(\left\{B \in \binom{[r]}{s} \mitte \bigcap_{i\in B} f(B\setminus \{i\})= \emptyset \right\}\right)>1-\varepsilon
		\end{equation*}
		follows.
	\end{theorem}
\end{remark}

\section{Reduction}
\label{Section_Reduction}
In this section, we construct the reduction from promise $(C_k,K_k)$ graph coloring to $\sigma$-gap $d$-multi-label cover with bound $B$.

\subsection{Direct Construction}
We give the concrete construction first and put it in the context of a Fourier transform in Section~\ref{Section_red_Fourier}.

For the reduction, we need to associate to every label cover instance a graph.
\begin{algdef}
	\label{algdef_reduction1}
	For a label cover instance 
	\begin{equation*}
		\inst =(n,m,C,(U_1,\dots,U_n),(U_{n+1},\dots,U_{n+m}), (i_c,j_c,\pi_c)_{c\in [C]} )
	\end{equation*}
	and an odd positve integer $k$ construct a graph as follows:
	
	Recall that powers of $C_k$ are defined by iterated products. Now, take the quotient of the graph $\bigcup_{i=1}^{n+m} C_k^{U_i}$ where for each $c\in [C]$ every tuple $(v_i)_{i\in U_{j_c}} \in C_k^{U_{j_c}}$, is identified with $(v_{\pi_c(i)})_{i\in U_{i_c}} \in C_k^{U_{i_c}}$.
\end{algdef}

\begin{theorem}
	\label{TheoremReduction1}
	Let $\inst $ be a label cover instance with bound $B$ and $k\in \IN, k\ge 3$ odd. Then, Definition~\ref{algdef_reduction1} constructs a graph in polynomial time. Moreover,
	\begin{enumerate}
		\item if $\inst $ is satisfiable, then $G$ maps homomorphically to $C_k$ and
		\item there are constants $\bar \sigma\in (0,1)$ and $\bar A\in \IN\setminus\{0\}$ (depending only on $k$ and not on $\inst$) such that if $\inst$ is a negative instance of the $\bar \sigma$-gap $\bar A$-multi-label cover problem, then there is no homomorphism from $G$ to $K_k$. 
	\end{enumerate}
\end{theorem}
We proof this theorem in the next main section. The constants will be chosen in Definition~\ref{algdef_constants}.

\subsection{Fourier Transform Perspective}
\label{Section_red_Fourier}
The reduction uses a construction similar to a finite Fourier transform or Pontryagin duality. Before explaining correctness of the reduction, we give this alternative perspective and properly define chosen constants.

\begin{definition}
	Consider a label cover instance 
	\begin{equation}
		\inst =(n,m,C,(U_1,\dots,U_n),(U_{n+1},\dots,U_{n+m}), (i_c,j_c,\pi_c)_{c\in [C]} )
		\label{eq_labelCover}
	\end{equation}
	with bound $B$.
	
	The \de{continuous dual} $\cdual \inst $ of $\inst $ is the infinite label cover instance
	$$
	(m,n,C,
	(\cdual U_{n+1},\dots,\cdual U_{n+m}),
	(\cdual U_1,\dots,\cdual U_n),
	(j_c,i_c,\cdual \pi_c)_{c\in [C]} 
	)
	$$
	where $\cdual U_i$ is $\map(U_i,\IT)$ for all $i\in [n+m]$ and $\cdual \pi_c$ is the map $\cdual U_{j_c} \to \cdual U_{i_c}$ defined by precomposition as $f\mapsto f\circ \pi_c$.
	
	The \de{finite dual} $\fdual \inst $ of $\inst $ with positive integer $D$ as step number is the label cover instance with bound $D^B$ given by
	$$
	(m,n,C,
	(\fdual U_{n+1},\dots,\fdual U_{n+m}),
	(\fdual U_1,\dots,\fdual U_n),
	(j_c,i_c,\fdual \pi_c)_{c\in [C]} 
	)
	$$
	where $\fdual U_i$ is $\map(U_i,\IT_D)$ for all $i\in [n+m]$ and $\fdual \pi_c$ is the map $\fdual U_{j_c} \to \fdual U_{i_c}$ defined by precomposition as $f\mapsto f\circ \pi_c$.
	
	The double dual of $\inst $ is defined as $\fdual{\fdual {\inst }}$.
\end{definition}
\begin{remark}
Note that the finite dual depends on the step number $D$ but we omit it in the notation. Note moreover:
\begin{enumerate}
	\item The finite dual with constant $D$ has indeed bound $D^B$.
	\item The finite dual is a subset of the continuous dual in the sense that $\fdual U_i\subseteq \cdual U_i$ holds for all $i$ and $\cdual \pi_c$ restricted to $\fdual U_{j_c}$ equals $\fdual \pi_c$. Every solution of $\fdual \inst $ is also a solution of $\cdual \inst $.
	\item The duals $\fdual U_i$ respectively $\cdual U_i$ can be considered as products of $\IT$.
	\item The finite dual $\fdual \inst $ has for all $x\in \IT_D$ the constant solution consisting of constant maps to $x$.
	\item In general, the double dual of $\inst $ will usually be different from $\inst $ and always solvable by the previous item.
\end{enumerate}
\end{remark}

\begin{definition}
	The duals $\fdual U_i$ respectively $\cdual U_i$ are metric spaces and groups where the group structure induced by the product structure from $\IT$ respectively $\IT_D$ and the metric is the supremum metric on the product. 
	
	Moreover, define for all $i\in [n]$ the map $T\colon \cdual U_i \to \cdual U_i$ to be the addition of $1/2$ in every component. Similarly define for $D$ even also $T\colon \fdual U_i \to \fdual U_i$. These maps are distinguished by the input.
\end{definition}
\begin{remark}
	\
	\begin{enumerate}
		\item The metric and the group structure are well-defined.
		\item The map $T$ is an isometry and satisfies $T^2=\id$. 
		\item We do not require the dual construction to be a homomorphism of any kind.
	\end{enumerate}
\end{remark}
Solutions of the dual or double dual are not necessarily related to solutions of the original instance.

For the reduction, we need to associate to every label cover instance a graph.
\begin{algdef}
	\label{algdef_reduction2}
	For a label cover instance $\inst $, a positive real number $q<1/2$ and a positive even integer $D$ construct a graph with a pseudo-metric on the vertices:
	Denote the label cover instance as in \eqref{eq_labelCover}.
	
	Set $V(G)$ to be the set $\bigcup_{i=1}^n \fdual U_i\cup \bigcup_{j=1}^m \fdual V_j$.
	
	Define the direct distance of between two points $v_1,v_2\in V(G)$ by the supremum metric in $U_i$ respectively $V_j$ if both points are in the same set. If they are in different sets, define it to be $0$, if there is $c\in [C]$ such that $\fdual\pi_c(v_1)=v_2$ or $\fdual\pi_c(v_2)=v_1$ and $\infty$ else. Now, define $\dist\colon V(G)^2 \to [0,\infty]$ to be the shortest path pseudo-metric.
	
	Define $T\colon V(G)\to V(G)$ to be the extension of the maps in $U_i$ or $V_j$.
	
	Set $E(G)$ to be 
	$
	\{(v_1,v_2)\mid \dist(v_1,T(v_2))\le q\}
	$.
\end{algdef}
For fixed $q\in \IQ\cap (0,1/2)$ and $D$ even, this algorithm computes $(V(G),E(G))$ in polynomial time in the size of $\inst $.

\begin{lemma}
	\label{Lemma_Reduction2}
	Let $\inst $ be a label cover instance and $k\in \IN, k\ge 3$ odd. Assume that $D$ is a multiple of $2k$ and that $q=1/(2k)$. Then,
	the graphs defined in Definition~\ref{algdef_reduction1} and~\ref{algdef_reduction2} are homomorphically equivalent.
\end{lemma}
\begin{proof}
	In fact, the graph structure on $\fdual U_i$ is $K_{\frac{D}{q(k-1)D}}^{U_i}$ in this case where $q(k-1)D$ is an integer. Note that if one would replace $\fdual U_i$ by $C_k^{U_i}$ the graph from Definition~\ref{algdef_reduction2} would be exactly the graph from Definition~\ref{algdef_reduction1}. (With the slight difference that for $\pi\colon U_i\to U_j$ and $x\in U_j$ the points $x$ and $\fdual\pi(x)$ are identified in Definition~\ref{algdef_reduction1} and just get the same neighbors in Definition~\ref{algdef_reduction2} which is irrelevant up to homomorphic equivalence.) Now, the claim follows from the homomorphic equivalence of
	$K_{\frac{D}{q(k-1)D}}$, $K_{\frac{2k}{k-1}}$ and $C_k$.
\end{proof}

We also give another generalization:
\begin{definition}
	We define the product with $r$-layers of a label cover instance $\inst $ as in \eqref{eq_labelCover} to be the $r-1$ fold product
	$$
	(U_{l_1},\dots,U_{l_{r-1}})_{l_1,\dots,l_{r-1} \in [n+m]}
	$$
	of $(n+m)^{r-1}$ sets with maps between them whenever all but one coordinates are equal and in the last coordinate, some $\pi_c$ is applied. For $i\in [r]$, the \de{principal layer} $i$ consists of tuples $(U_{l_1},\dots,U_{l_r})$ where $l_1,\dots,l_{i-1}\le n$ and $l_{i},\dots,l_{r-1}>n$.
	
	The continuous dual and the finite dual of the product with layers are defined analogously to the duals of label cover instances.
	
	For a label cover instance $\inst $ and positive integers $k$ (odd), $r\ge 2$ and $D$ (a multiple of 2k), define $\ix_{k,D,r}$ to be the graph obtained from the finite dual of the product with $r$ layers of $\inst $ as in Definition~\ref{algdef_reduction2} where $q=1/(2k)$. Additionally, define $\iix_{k,r}$ to be the set obtained from the continuous dual. Define $T\colon \iix_{k,r} \to \iix_{k,r}$ similar to the map $T$ from above.
\end{definition}

\begin{lemma} \label{lemma_homequiv}
	For a label cover instance $\inst $ and positive integers $k$ (odd), $r\ge 2$ and $D$ (a multiple of 2k), the graphs $G$ defined in Definition~\ref{algdef_reduction1} and $\ix_{k,D,r}$ are homomorphically equivalent.
\end{lemma}
\begin{proof}
	By Lemma~\ref{Lemma_Reduction2}, $G$ and $\ix_{k,D,2}$ are homomorphically equivalent. The homomorphic equivalence of $\ix_{k,D,r}$ and $\ix_{k,D,2}$ is induced by the projection 
	$$
	(U_{l_1},\dots,U_{l_{r-1}})_{l_1,\dots,l_{r-1} \in [n+m]} \to U_{l_1}
	$$ 
	respectively the diagonal embedding $U_{l} \to (U_{l},\dots,U_{l})$.
\end{proof}

To simplify the upcoming parts of the proof, they use different variables.
\begin{algdef} \label{algdef_constants}
	Given a positive odd integer $k\ge 3$, define additional variables in the following way:
	\begin{enumerate}
		\item Define $\bar q=1/(2k)$.
		\item Choose $\bar \varepsilon\in \IQ$ with $0<\bar \varepsilon<1/(8k)$.
		\item Choose $\bar L\in \IN$ with $\bar L \ge 4k$.
		\item Define $\bar \beta\in \IQ$ and $\bar K\in \IN$ by Theorem~\ref{TheoremJunta3} respectively
		Remark~\ref{RemarkJuntaQuantitative} depending on $L=\bar L$ and $\varepsilon=\bar\varepsilon$.
		\item Choose $\bar m\in \IN$ with $\bar m>4\bar L/\bar \beta$.
		\item Choose $\bar \lambda\in \IQ$ with $0<\bar \lambda<1, \bar \lambda<\bar \varepsilon/\bar K$.
		\item Choose $\bar s\in \IN$ with 
		$(1-\frac{\bar \lambda^{\bar m}}{\sum_{l=0}^{\bar m}\bar \lambda^l})^{\bar s}<\bar \varepsilon, \bar s\ge 2$ 
		\item Choose $\bar \gamma\in \IQ$ with $\bar \gamma>0, \bar \gamma^2<\bar \varepsilon\bar \beta^2/(4s)$.
		\item Choose $\bar P\in \IN$ with $\bar P$ divisible by $2k(k-1)$, 
		$\bar P>2\bar L/(\bar \gamma \bar m),\bar P > 2\bar L\bar s/\bar \varepsilon$.
		\item Define $\bar D \in \IN$ as $\bar m\bar P$.
		\item Define $\bar A\in \IN$ as $J(\bar s\bar L,\bar \gamma)(\bar m+1)^{\bar s}\ge 1$ where $J$ is the Junta function from Theorem~\ref{TheoremJunta2}. 
		\item Define $\bar r\in \IN$ and $\mu$ as probabillity measure on $\bar s$ element subsets of $[\bar r]$ according to Lemma~\ref{Lemma_Ramsey} with $s=\bar s, d=\bar A, \varepsilon=\bar \varepsilon$.
		\item Choose $\bar \sigma\in \IQ$ with  $0<\bar \sigma<1$ and $\bar \sigma<\bar \varepsilon/(4^{\bar r})$.
		\item Choose $\bar B$ such that the $\bar \sigma$-gap $\bar A$ multi-label cover problem with bound $\bar B$ is NP hard (see Corollary~\ref{Corollary_BoundLabelCover}/Remark~\ref{Remark_BoundLabelCoverQuantitative}).
	\end{enumerate}
\end{algdef}
None of these variables depends on any label cover instance.
\begin{remark}
	The variables are denoted similar to \cite{OAI_PCSP}, a bar is added to mark this set of variables and some estimations are simplified.
	Additionally, one restriction to $\bar \beta$ and $\bar L$ is added and needed for the generalization from $K_3$ to odd cycles. 
\end{remark}

\section{Correctness}
\label{Section_Correctness}
This section is dedicated to prove Theorem~\ref{TheoremReduction1} and concluding Theorem~\ref{TheoremMain}.
\subsection{Satisfiable Instance}
\label{Section_corr_SatInst}
The easy direction.
\begin{lemma}\label{lemma_easyDirection}
	Let $\inst $ be a satisfiable label cover instance and $k\in \IN, k\ge 3$ odd. Then, Definition~\ref{algdef_reduction1} (where $q=1/(2k)$ and $D$ is a multiple of $2k$), defines a graph $G$ that maps homomorphically to $C_k$.
\end{lemma}
\begin{proof}
	Choose a solution of $\inst $. The respective projections induce maps $C_k^{U_i}\to C_k$. These maps map the identified points in $G$ to the same point in $C_k$ as the solution of $\inst $ enforces compatibility with the $\pi_c$ maps.
	
	In categorical language, the 1-element set $S$ maps homomorphically into $\lim \inst $ by assumption. By applying $\Hom(\blank, C_k)$, we get that there is a graph homomorphism from $G=\colim_{U\in \inst } \Hom(U,C_k)$ to $C_k=\Hom(S,C_k)$.
\end{proof}

\subsection{Probability Distribution}
\label{Section_corr_ProbDist}
The proof for the unsatisfiable instances is much harder as already teased in the overview in Section~\ref{Section_overview}.
Throughout the following sections, we fix $k$, the constants from Definition~\ref{algdef_constants} and a label cover instance $\inst $ with bound $\bar B$ denoted as in \eqref{eq_labelCover}.

We define $\ix$ to be the graph $\ix_{k,\bar D,\bar r}$ where $\bar D$ and $\bar r$ are chosen according to Definition~\ref{algdef_constants}.
If $\ix$ has a loop, then there is clearly no graph homomorphism $\ix\to C_k$. For the following sections we assume that no loop exists.

Some elements related to the layered product and its dual get a new shorter name.
\begin{definition}
	Define $W_-$ to be the product $[C]^{\bar r-1}\times \binom{[\bar r]}{\bar s}$.
	
	For any element $w=(H,B)$ in $[C]^{\bar r-1}$ and $a\in [\bar r]$, define 
	$$
	M_{w,a}=M_{H,a}=(U_{j_{c_1}},\dots,U_{j_{c_{a-1}}},U_{i_{c_{a}}},\dots,U_{i_{c_{\bar r-1}}} )
	$$
	as the set in the $\bar r$-layered product of $\inst $, in the main $a$-th layer. For $w=(c_1,\dots,c_{r-1},B)\in W_-$, define $M_w=M_{(c_1,\dots,c_{r-1}),\min(B)}$.
	
	For any element $H$ in $[C]^{r-1}$ and $a,b\in [r]$, $a<b$, define $\pi_{(H,a,b)}$ (or $\pi_{(w,a,b)}$ if the element comes from $w\in W_-$) to be the map in the $r$ layered product of $\inst $ from $M_{w,a}$ to $M_{w,b}$ that composes by $\pi_{c_i}$ for all $a\le i<b$ and by identity in the other components. By $\fdual M_w$ and $\cdual M_w$ we denote the spaces $\IT_D^{M_w}$ respectively $\IT^{M_w}$ which we consider as subsets of the dual of the layered product of $\inst $ in the natural way and thus as subsets of the vertices of $\ix$.
\end{definition}

We define some probability measures which add up to a measure on the vertices of $\ix$.
\begin{algdef}\label{algdef_measures}
	Some probability measures on finite sets.
	\begin{enumerate}
		\item Consider the counting measure on $[C]^{\bar r-1}$. Divided by $C^{\bar r-1}$, this imposes a probability measure $\mu_H$.
		\item Consider the measure $\mu$ on $\binom{[\bar r]}{\bar s}$ from Lemma~\ref{Lemma_Ramsey} that was already defined in Definition~\ref{algdef_constants}.
		\item Define $\mu_W$ to be the probability distribution on $W_-$ that chooses a tuple in $[C]^{\bar r-1}$ according to the equal distributed measure and $B\in \binom{[\bar r]}{\bar s}$ according to $\mu$.
		\item Let $\rho$ be the probability distribution on $\{0,1/\bar m,2/\bar m,\dots, \bar m/\bar m\}$ where the probability of $i/\bar m$ is given by $p_i=\frac{\lambda^i}{\sum_{l=0}^{\bar m}} \lambda^l$.
		\item Let $\rho_\theta$ be the probability distribution on $\{0,1/\bar P,2/\bar P,\dots, (\bar P-1)/\bar P\}$ that chooses all results equally by probability $1/\bar P$.
		\item Denote by $\lebesque$ the standard Lebesque measure on $\IR$ and the probability measure obtained on $\IT=\IR/\IZ$ by first taking the restriction of the measure to $[0,1)$ and then taking the quotient.
		\item For every $w=(H,B)\in W_-$, $I\subseteq B$ non-empty, define $x_{H,I}$ as
		\begin{align*}
			[0,1]^I \times [0,1]^{\prod_{j\in I} M_{w,j}} &\to \IT^{M_{w,\min I}} \\
			(t,\theta) &\mapsto x_{H,I}(t,\theta)
			\\
			(t_j, (\theta_{j,l})_{l\in M_{w,j}})_{j\in I} 
			&\mapsto \left(\left(\sum_{j\in I} t_j \theta_{j, \pi_{(w,\min I,j)}(\kappa)}\right)\!\!\!\pmod {1}\right)_{\kappa \in M_{w,\min I}}
		\end{align*}
		and also denote its restriction 
		$$
			([0,1]\cap 1/\bar m\IZ)^I \times ([0,1)\cap 1/\bar P\IZ)^{\prod_{j\in I} M_{w,j}} \to \IT_D^{M_{w,\min I}}
		$$
		the same way. Define $x_w=x_{H,B}$.
		Now, define $\nu_w$ as the probability measure on $\fdual M_{w}\subseteq V(\ix)$ defined by
		$$
			\forall x \in \fdual M_w\colon 
			\nu_w(\{x\})= \rho^B \times \rho_\theta^{\prod_{j\in B} M_{w,j}} (\{(t,\theta) \mid x=x_w(t,\theta) \}) 
		$$
		on atoms.
		\item Finally, define $\nu$ to be the probability distribution on $V(\ix)$ as $$\nu(X)=\sum_{w\in W_-}(\nu_w(X) \cdot \mu_W(\{w\})).$$  
	\end{enumerate}
\end{algdef}
\begin{remark}
	The interesting part of this definition is the definition of $x_w$. 
	The operation $\pmod{1}$ denotes the quotient map $\IR\to\IT$. 
	Considering $t$ as a parameter we consider the image of $x_w$. The choice $t=0$ assigns 0 to every input $\theta$. 
	When $t_l=1$ (and all other entries stay 0), then the image becomes instead all points $(x_1,\dots,x_{|M_{w,\min B}|})$ in $\IT^{M_{w,\min B}}$ where $x_i=x_j$ whenever $i$ and $j$ have equal images under $\pi_{(w,\min B,l)}$. 
	Having $t\in (0,1)$ induces as image only points close to zero and having multiple nonzero indices allows for the sum the images.
\end{remark}

\subsection{Independent Sets}
\label{Section_corr_IndependentSet}
In the last section, we have fixed many constants and defined new probability measures.
For the rest of the correctness proof fix any independent set $\mathcal A$ in $V(\ix)$. We show that if $\inst $ is a negative $\bar \sigma$-gap $\bar A$-multi-label cover instance, then $\nu(\mathcal A)<1/k$ by doing multiple analytic approximations. It implies that there is no graph homomorphism $\ix \to K_k$ as the preimage of the $k$ vertices are independent sets of which at least one has measure $\ge 1/k$.

\begin{remark}\label{remark_cliqueResult}
	In fact, this even implies that we can define a graph $\ix'$ that has no clique of size $|V(\ix')|/k$. We will discuss this further in the appendix.
	
	Construct the graph $\ix'$ from $\ix$ by replacing every vertex $v$ of $\ix$ with $\nu(v)\cdot g$ many vertices, where $g\in \IN$ is chosen as the smallest positive integer that makes $\nu(v)\cdot g$ an integer for all $v$. Recall that every measure of a point is rational and the denominator can be bounded independently of the instance. Two vertices in $\ix'$ are connected if their original vertices in $\ix$ are. In fact, for fixed $k$ and assuming that $\inst$ is bounded by $\bar B$, this graph $\ix'$ can be constructed in polynomial time from the label cover instance (assuming that every set in the instance is bounded).
\end{remark}

To estimate $\nu(\mathcal A)$, we define a Lipschitz-continuous function that is similar to a characteristic function.
We start by defining a metric.
\begin{definition}
Recall that $\fdual \inst $ and similar $V(\ix)$ carries a quasi-metric. We extend this metric to $\iix\coloneqq \iix_{k,\bar r}\supseteq V(\ix)$ by choosing as local distance function $\distloc(x,y)=\max_{i\in M}(\dist(x_i,y_i))$ if $x$ and $y$ are from the same set $\IT^M$ and $\infty$ otherwise and we define the metric globally as $\dist(x,y)\coloneqq $
\begin{align*}
	\min(&\{\distloc(x,y)\}\cup{}
	\\
	&\{\distloc(x,x_1)+\dist_{V(\ix)}(x_1,y_1)+\distloc(y_1,y)\mid x_1,y_1\in V(\ix)\}).
\end{align*}
\end{definition}
Note that this definition actually defines a quasi-metric.
\begin{remark} \label{RemarkSmallDistanceToPi}
	While for all maps $\fdual\pi$ induced by maps in $\inst $ the condition implies in $\ix$ that $\dist(x,\fdual\pi(x))=0$, we get in this metric in $\iix$ only that 
	\begin{align*}
		\dist(x,\cdual\pi(x)) \le 1/ \bar D
	\end{align*}
	holds.
	The advantage is that the distance of two points in $V(\ix)\subseteq \iix$ equals the distance in $V(\ix)$.
\end{remark}
\begin{proof}
	A point $x$ in $\iix$ has always distance at most $1/(2 \bar D)$ to a point $x_1$ in $V(\ix)$ that has distance $0$ to its $\cdual\pi$-image $\cdual\pi(x_1)=\fdual\pi(x_1)$ that again differs only by $1/(2 \bar D)$ from $\cdual\pi(x)$.
\end{proof}
\begin{definition}
	\label{Definition_mapF}
	Define the map $f=f_\mathcal{A}$ as
	\begin{align*}
		f\colon \iix &\to [-1,1] \\
		x &\mapsto 
		\begin{cases}
			1-\bar L\dist(x,\mathcal{A})&\text{ if }	1-\bar L\dist(x,\mathcal{A})>0\\
			\bar L\dist(Tx,\mathcal{A})-1&\text{ if }	1-\bar L\dist(Tx,\mathcal{A})>0 \\
			0&\text{ else}
		\end{cases}
	\end{align*}
\end{definition}
\begin{lemma}
	For every points $y_1,y_2\in \mathcal{A}$, $$\dist(y_1,T(y_2))>1/(2k)$$ holds. Moreover, the map $f$ in Definition~\ref{Definition_mapF} is well-defined and $L$-Lipschitz continuous.
\end{lemma}
\begin{proof}
	By the definition of edges in $\ix$, there is no edge between $y_1$ and $y_2$ if and only if $\dist(y_1,T(y_2))>\bar q = 1/(2k)$. As there are no edges in $\mathcal{A}$, the result follows.
	For the definition, note that if both $1-\bar L\dist(x,\mathcal{A})>0$ and $1-\bar L\dist(Tx,\mathcal{A})>0$ hold, then 
	$$
		1/(2k)<\dist(\mathcal{A}, T(\mathcal{A}))\le \dist(\mathcal{A},x)+\dist(\mathcal{A}, Tx)<2/\bar L\le 1/(2k)
	$$
	where the last inequality follows from $\bar L \ge 4k$. As this is a contradiction, there is no such intersection. The Lipschitz continuity follows similar.
\end{proof}
\begin{remark}
	The following holds.
	\begin{align*}
		\forall x\in \mathcal{A}  &\colon & f(x)&=1
		\\
		\forall x\in \iix &\colon& f(x)&=-f(T(x))
		\\
		\forall x\in \mathcal{A}&\colon & f(T(x))&=-1
	\end{align*}
\end{remark}

\begin{remark} \label{Remark_BaseTransition}
	For $x\in \IT^{M_{w,1}} \subseteq \iix$ and $\cdual\pi\colon \IT^{M_{w,1}} \to \IT^{M_{w,2}}$ a map induced by a map of the label cover instance, the inequality
	$$
		|f(x)-f(\cdual\pi(x))|\le \bar L \dist(x,\cdual\pi(x))<\bar L / \bar D = \bar L / (\bar P \bar m)<\bar \gamma/2
	$$
	follows using Remark~\ref{RemarkSmallDistanceToPi}.
\end{remark}

\subsection{Junta Approximation}
\label{Section_corr_JuntaApprox}
Now, we approximate the map $f(x(t,\theta))$ by maps that depend only on few coordinates.
\begin{definition}
	For fixed $w=(H,B)\in W_-$, $I\subset B$ non-empty, $t\in [0,1]^I$, define
	\begin{align*}
		h_{H,I,t}\colon \prod_{j\in I} [0,1]^{M_{w,j}} &\to [-1,1]\\
		\theta &\mapsto f(x_{H,I}(t,\theta))
	\end{align*}
\end{definition}
Note that $h_{H,I,t}$ does not depend on $B$.
\begin{lemma}[vgl. Lemma 4.1 in \cite{OAI_PCSP}]
	\label{Lemma_ManyJuntaApproximations}
	For every $w=(H,B)\in W_-$, $I\subset B$ nonempty, there are sets $S_{H,I,j}\subset M_{w,j}$ for all $j\in I$ with
	$$
		\sum_{j\in I} |S_{H,I,j}|\le \bar A=J(\bar s\bar L,\bar \gamma)(\bar m+1)^{\bar s}
	$$
	and maps 
	\begin{align*}
		g_{H,I,t}\colon \prod_{j\in I} [0,1]^{M_{w,j}} &\to [-1,1]
	\end{align*}
	that only depends on coordinates $(j,l)$ with $l\in S_{H,I,j}\subset M_{w,j}$ and satisfy
	$$
		\int_{\prod_{j\in I} [0,1]^{M_{w,j}}} (h_{H,I,t}(\theta)-g_{H,I,t}(\theta))^2\deld\lebesque^{\coprod_{j\in I} M_{w,j}}(\theta)<\bar \gamma^2
	$$
	for all $t\in ([0,1]\cap 1/\bar m\IZ)^I$.
\end{lemma}
\begin{proof}
	Apply Theorem~\ref{TheoremJunta2} to every of the $(m+1)^s$ possible choices of $t$.
\end{proof}
To make the sets $ S_{H,I,j}\subset M_{w,j}$ better comparable for different $I$, we map them into $M_{w,r}$.
\begin{definition}
	For $H\in [C]^{r-1}, I\in \powerset^+_{\le \bar s}([\bar r])$, define $A_{H,I} \subseteq M_{H,r}$ as the union $\cup_{j\in I}\pi_{(H,j,r)}(S_{H,I,j})$.
\end{definition}

\subsection{Applying Negative Instances}
\label{Section_corr_NegInst}
Up next, we use the assumption of a negative instance to estimate measures.

\begin{lemma}[Lemma 4.2 of \cite{OAI_PCSP}] \label{Lemma_AssumptionNegInst}
	Assume that $\inst$ is a negative instance of the $\bar \sigma$-gap $\bar A$-multi-label cover problem with bound $B$. Choose $I,J\in \powerset^+_{\le \bar s}([\bar r])$ with $\max(I)<\min(J)$. Then, the bound
	$$
		\mu_H(\{H\in [C]^{\bar r -1}\mid A_{H,I} \cap A_{H,J} \ne  \emptyset\})\le \bar \sigma
	$$
	holds.
\end{lemma}
\begin{proof}
	Apply the bound assumption in the right way.
	
	We fix a tuple $H'=(c'_1,\dots,c'_{\bar r -1})\in [C]^{\bar r -1}$ and show that 
	\begin{align*}
	|\{H=(c_1,\dots,c_{\bar r -1})\in [C]^{\bar r -1}&\mid A_{H,I} \cap A_{H,J} \ne  \emptyset, \\&\phantom{{}\mid{}}
	\forall j<\bar r\colon j= \min J - 1 \lor c_j=c'_j
	\}|\le \bar \sigma \cdot |C|
	\end{align*}
	holds. By summing over all $H'$, the result follows.
	
	Now, consider the label cover instance $\inst$. We construct for every $U_i$ a subset $U'_i\subseteq U_i$ in the following way: 
	Recall the definitions of $S_{H,I,a}$ and $S_{H,J,a}$ with $a\in I$ respectively $a\in J$. Now, we can project $S_{H,I,a}\subseteq M_{H,a}$ to $U_{i_{c_{(\min J-1)}}}$ by taking the $(\min J-1)$-th projection. Similarly, $S_{H,J,a}$ projects to $U_{j_{c_{(\min J-1)}}}$.
	$$
	\begin{tikzcd}
		S_{H,I,a} 
			\ar[r,hook] 
			\ar[d] 
		&M_{H,a}=(U_{j_{c_1}},\dots,U_{j_{c_{a-1}}},U_{i_{c_{a}}},\dots,U_{i_{c_{\bar r-1}}} )
			\ar[d, "\pi_{(w,i,r)}"']
			\ar[r]
		&U_{i_{c_{(\min J-1)}}}
			\ar[d, "\pi_{c_{(\min J-1)}}"']
		\\
		A_{H,I}
			\ar[r,hook] 
		&M_{H,r}=(U_{j_{c_1}},\dots,U_{j_{c_{\bar r-1}}})
			\ar[r]
		& U_{j_{c_{(\min J-1)}}}
		\\
		%second part
		S_{H,J,a} 
			\ar[r,hook] 
			\ar[d] 
		&M_{H,a}=(U_{j_{c_1}},\dots,U_{j_{c_{a-1}}},U_{i_{c_{a}}},\dots,U_{i_{c_{\bar r-1}}} )
			\ar[d, "\pi_{(w,i,r)}"']
			\ar[r]
		&U_{j_{c_{(\min J-1)}}}
			\ar[d, "="']
		\\
		A_{H,J}
			\ar[r,hook] 
		&M_{H,r}=(U_{j_{c_1}},\dots,U_{j_{c_{\bar r-1}}})
			\ar[r]
		& U_{j_{c_{(\min J-1)}}}
	\end{tikzcd}
	$$
	We define $U'_{H,i_{c_{(\min J -1)}}} \subseteq U_{i_{c_{(\min J-1)}}}$ to be the union of these sets $\cup_{a\in I} S_{H,I,a}$ and $U'_{H,j_{c_{(\min J -1)}}} \subseteq U_{j_{c_{(\min J-1)}}}$ to be the union $\cup_{a\in J} S_{H,J,a}$. Note that we can choose $S_{H,I,a}$ independent of $c_{j}$:
	\begin{itemize}
		\item If $H$ and $H''$ agree on every component, except for $c_{\min J-1}$ but $i_{c_{\min J-1}} = i_{c''_{\min J-1}}$ then we may assume $S_{H,I,a}=S_{H'',I,a}$ for all $a\in I$ and thus get that $U'_{H,i_{c_{(\min J -1)}}}$ is independent of $H$.
		\item If $H$ and $H''$ agree on every component, except for $c_{\min J-1}$ but $j_{c_{\min J-1}} = j_{c''_{\min J-1}}$ then we may assume $S_{H,J,a}=S_{H'',J,a}$ for all $a \in J$ and thus get that $U'_{H,j_{c_{(\min J -1)}}}$ is independent of $H$.
		\item We have $i_{c_{(\min J-1)}}\le n$ and $j_{c_{(\min J-1)}}>n$ thus the two cases do not interfere.
	\end{itemize} 
	Now, we defined for every set $U_i$ of the instance a subset $U'_i=U'_{H,i}$ of size at most $\bar A$ that is independent of $H$. (If there is no $H$ we note that this set is not contained in any constraint.)
	As $\inst$ is a negative instance, we get $\pi_c(U'_{i_c})\cap U'_{j_c} =\emptyset$ for all constraints but an $\sigma$-fraction.
	
	Note that $A_{H,I}\cap A_{H,J}\ne \emptyset$ implies that their images in $U_{j_{(\min J-1)}}$ also intersect nontrivial, which is exactly the statement 
	$$
		\pi_{c_{(\min J-1)}}(U'_{H,i_{\min J -1}}) \cap 
		U'_{H,j_{\min J -1}} \ne \emptyset
	$$
	which occurs only in an $\sigma$-fraction by the previous paragraph.
\end{proof}
\begin{corollary}
	\label{Corollary_NegInst}
	Assume that $\inst$ is a negative instance of the $\bar \sigma$-gap $\bar A$-multi-label cover problem with bound $\bar B$.
	Then,
	$$
	\mu_H(\{H\in [C]^{\bar r -1}\mid \forall I,J\colon\max(I)<\min(J) \implies  A_{H,I} \cap A_{H,J} \ne  \emptyset\})\le 4^{\bar r}\bar\sigma<\bar \varepsilon
	$$
	where $I$ and $J$ are taken from $ \powerset^+_{\le \bar s}([\bar r])$ as usual.
\end{corollary}
\begin{proof}
	There are at most $2^r$ possible sets $I$ and $J$ each as they are in $\powerset^+([\bar r])$. Summing Lemma~\ref{Lemma_AssumptionNegInst} over all possibilities gives the desired estimation. The final inequality follows from Definition~\ref{algdef_constants}.
\end{proof}

\subsection{Applying Ramsey Argument and Aligned Elements}
\label{Section_corr_Ramsey}
By design, the results from the previous section are exactly the assumptions needed in the Ramsey part and they allow us to estimate further and exclude everything but aligned pairs.
\begin{definition}
	An element $w=(H,B)\in W_-$ is \de{aligned} if for all $I,J\in \powerset^+(B)$, we have $A_{H,I}\cap A_{H,J}=\emptyset$.
	In this case, there is a map $\alpha\colon M \to B$ such that 
	$$
	\forall I \in \powerset^+(B)\colon \alpha(A_I)\subseteq I.
	$$
	Define $W_+=\{(H,B,\alpha)\mid (H,B)\in W_-\text{ aligned with map } \alpha\}$.
\end{definition}
Note that $\alpha$ indeed exists and may depend on $B$. We may consider $W_+$ as subset of $W_-$ and extend $\mu_W$ to a measure on $W_+$ (where we might have $\mu_W(W_+)<1$) by choosing any probability distribution on such maps $\alpha$.
\begin{lemma} \label{Lemma_notAlligned}
	We have
	$$
	\mu_W(\{(H,B)\mid (H,B)\text{ aligned}\})>1-2\bar \varepsilon.
	$$
	whenever $\inst$ is a negative instance.
\end{lemma}
\begin{proof}
	By Corollary~\ref{Corollary_NegInst}, the measure of well-behaved $H$ is at least $1-\bar \varepsilon$. By Lemma~\ref{Lemma_Ramsey}, applied with $s=\bar s, d=\bar A$ and $\varepsilon=\bar \varepsilon$, each of those $H$ admits with $\mu$-measure at least $1-\bar \varepsilon$ aligned $B$-s. Thus, the measure of aligned pairs is at least $(1-\bar \varepsilon)^2>1-2\bar\varepsilon$.
\end{proof}

\begin{remark}
	To get some intuition about aligned tuples, consider the diagrams
	$$
	\begin{tikzcd}
		B&M_{\bar r} \ar[l,"\alpha"'] & \dots \ar[l, "\pi_{(w,\bar r-1,\bar r)}"'] & M_2 \ar[l, "\pi_{(w,2,3)}"'] & M_1 \ar[l, "\pi_{(w,1,2)}"']
		\\
		\IT^B \ar[r, "(\blank)_\alpha"]
		&\IT^{M_{\bar r}} \ar[r] \ar[drrr,"f"']
		&\dots \ar[r]
		&\IT^{M_2} \ar[r, "\cdual\pi_{(w,1,2)}"] \ar[dr,"f"]
		&\IT^{M_1} \ar[d, "f"]
		\\
		&&&& {[-1,1]}
	\end{tikzcd}
	$$
	which are commutative up to an $\bar \gamma / 2$-distance. However, the map $f$ is approximated on every torus $\IT^{M_{i}}$ by a different map $g$ that only depends on few coordinates. Now alignment is the property that these coordinates are different for each map to $[-1,1]$, and even their images in $M_{\bar r}$ are different. This should be unlikely if $f$ is not close to a constant function. A fact we make precise in the next section.
\end{remark}

\subsection{Analytic Methods}
\label{Section_corr_Analysis}
We want to show with analytic methods that for an aligned $w\in W_+$, the measure
$
	\nu_w(\mathcal{A})
$ is small. This combines with Lemma~\ref{Lemma_notAlligned} to a full estimation of $\nu(\mathcal{A})$.

We start by defining two variations of $\nu_w$.
\begin{definition}
	For $w=(H,B)$, define $\fdual\nu_w$ to also be the probability measure on 
	$$
		M_w \times [0,1]^B \times [0,1]^{\prod_{j\in B} M_j}
	$$
	as
	$$
		\fdual\nu_w(X) = \rho^B \times \rho_\theta^{\prod_{j\in B} M_j}
		(\{(t,\theta)\mid (x_{H,I}(t,\theta),t,\theta) \in X\})
	$$
	and define $\cdual\nu_w$ on the same space by
	$$
		\cdual\nu_w(X) = \rho^B \times \lebesque^{\prod_{j\in B} M_j}
	(\{(t,\theta)\mid (x_{H,I}(t,\theta),t,\theta) \in X\})
	$$
\end{definition}
Now, $\nu_w$ is the projection of $\fdual\nu_w$ and $\cdual\nu_w$ differs from $\fdual\nu_w$ by choosing $\theta$ in a power of $[0,1]$ continuously instead of discrete.

\begin{lemma} \label{lemma_fullCoefficientNeeded}
	Consider $w\in W_+$. Then,
	$$
		\fdual\nu_w(\{(x,t,\theta)\mid \nexists j\in B\colon t_j=1\})<\bar \varepsilon
	$$
\end{lemma}
\begin{proof}
	The constant $\bar s$ was chosen large enough. (We have
	$$
		\fdual\nu_w(\{(x,t,\theta)\mid \nexists j\in B\colon t_j=1\})
		= \prod_{j\in B} (1-\frac{\bar \lambda^{\bar m}}{\sum_{l=0}^{\bar m}\bar \lambda^l})
		=(1-\frac{\bar \lambda^{\bar m}}{\sum_{l=0}^{\bar m}\bar \lambda^l})^{\bar s}<\bar \varepsilon
	$$ 
	for the left side.)
\end{proof}

\begin{definition}
	For $w\in W_+$, define
	\begin{align*}
		H_w\colon [0,1]^B \times [0,1]^{\prod_{j\in B} M_j} \times \IT^B &\to [-1,1] \\
		(t,\theta,z)&\mapsto f(x_w(t,\theta)+(z_{\alpha(\pi_{(w,\min B,r)}(\kappa))})_{\kappa \in M_{\min B}} )
	\end{align*}
	where $+$ denotes the addition in $\IT^{M_{\min B}}$.
\end{definition}
This map can be considered as generalization of $f(x_w(t,\theta))$ that adds an additional summand in all coordinates where the summand is equal in coordinates with the same image in $B$.

Up next, we want to bound the dependence of $H_w$ on single $z$-coordinates.
\begin{definition}
	For $i\in B$ and $w\in W_+$, define 
	\begin{align*}
		\IE_{z_i}H_w(t,\theta,z)&\coloneqq \int_{\IT} H_w(t,\theta,z+e_i(z')) \deld \lebesque (z')
		\\
		\IE_{z}(H_w(t,\theta,z))^2&\coloneqq \int_{\IT^B} (H_w(t,\theta,z'))^2 \deld \lebesque^B (z')
	\end{align*}
	where $e_i(z)\in \IT^B$ is the vector with $z'$ in coordinate $i$ and $0$ everywhere else.
	
	Define 
	$$D_i(t,\theta)\coloneqq \sqrt{\int_{\IT^B} (H_w(t,\theta,z) - \IE_{z_i}H_w(t,\theta,z))^2\deld \lebesque^{B}(z)}$$
	as an $L^2$-norm if the influence of the $z_i$-component.
\end{definition}

\begin{lemma} \label{Lemma_SecondJuntaApplied}
	Consider $w\in W_+$. Now, for all $t\in [0,1]^B$ and all $\theta\in \prod_{j\in I} [0,1]^{M_j}$, it follows that
	$$
		\int_{\IT^B} H_w(t,\theta,z) \deld\lebesque^B(z)=0
	$$
	and moreover 
	\begin{align*}
		|\{i\in B\mid D_i(t,\theta)\ge \bar \beta/2\}|&\le \bar K,
		\\
		\IE_zH_w(t,\theta,z)^2>\bar\varepsilon \implies |\{i\in B\mid D_i(t,\theta)\ge \bar \beta\}|&\ge 1
	\end{align*}
\end{lemma}
\begin{proof}
	The first claim follows from $H_w(t,\theta,z) =-H_w(t,\theta,T(z))$ where $T$ adds $1/2$ in every component. The second claim follows from applying Theorem~\ref{TheoremJunta3} for $L=\bar L, \varepsilon=\bar \varepsilon, \beta=\bar\beta, K = \bar K, s=\bar s, H(\blank)=H_w(t,\theta,\blank)$.
\end{proof}
\begin{lemma} \label{Lemma_UsingAllinged}
	Consider $w\in W_+$. Now,
	$$ 
		\cdual\nu_w \{(x_w(t,\theta),t,\theta) 
		\mid \exists i,j\in B\colon t_j=1 \land t_i=0 \land D_i(t,\theta)\ge \beta \}\le \frac{4\bar \gamma^2 \bar s}{\bar \beta ^2}<\bar \varepsilon
	$$
	holds.
\end{lemma}
\begin{proof}
	Essentially using the aligned property.
	
	The right inequality follows from the choice of $\bar \gamma$. For the left inequality, note that we can define for each $z\in \IT^B$ a vector $z_{\alpha\circ \pi_{(w,\min B,r)}}\in \IT^{M_{\min B}}$ as a minor. We consider $\IT^{M_{\min B}}$ as a subset of $[0,1]^{\prod_{j\in B} M_j}$ by taking the preimage $[0,1)\to \IT$ and $0$ in other coordinates and get
	$$
		H_w(t,\theta,z)=H_w(t,\theta+z_{\alpha\circ \pi_{(w,\min B,r)}}\pmod{1},0)
	$$
	by unfolding the definition. Note that since $t_j=1$, this map is for this fixed $t$ continuous in $z$ even though the $(\blank)\pmod{1}$-operation is not. We omit writing the $(\blank)\pmod{1}$-operation in this situation in the rest of the proof.
	We have
	$$
		H_w(t,\theta+z_{\alpha\circ \pi_{(w,\min B,r)}},0)
		=f(x_{H,B}(t,\theta+z_{\alpha\circ \pi_{(w,\min B,r)}}))
	$$
	and with $t'$, $\theta'$ and $z'$ denoting the restrictions of $t$, $\theta$ and $z_{\alpha\circ \pi_{(w,\min B,r)}}$ to $B\setminus\{i\}$, we also have 
	\begin{align*}
			f(x_{H,B}(t,\theta+z_{\alpha\circ \pi_{(w,\min B,r)}}))
		=f(\cdual\pi_{(w,\min B\setminus\{i\}, \min B)}(x_{H,B\setminus\{i\}}(t',\theta'+z')))
	\end{align*}
	where $t'$, $\theta'$ and $z'$ denote the restrictions of $t$, $\theta$ and $z_{\alpha\circ \pi_{(w,\min B,r)}}$ to $B\setminus\{i\}$. We may now invoke Remark~\ref{Remark_BaseTransition} to get that
	$$
		|f(x_{H,B}(t,\theta+z_{\alpha\circ \pi_{(w,\min B,r)}}))
		-f((x_{H,B\setminus\{i\}}(t',\theta'+z')))
		|
		\le \bar\gamma/2<\bar \gamma
	$$
	holds. Moreover, by choosing $I=B\setminus\{i\}$ we get for all relevant $t$, that are $t\in ([0,1]\cap 1/\bar m\IZ)^I$ with $t_j=1$, that we can apply Lemma~\ref{Lemma_ManyJuntaApproximations} and get
	$$
		\int_{\prod_{j\in I} [0,1]^{M_j}}
			(
				f((x_{H,B\setminus\{i\}}(t',\theta')))
				-g_{H,I,t'}(\theta')
			)^2 
		\deld\lebesque^{\coprod_{j\in I} M_j}(\theta')<\bar \gamma^2
	$$
	which implies, as adding $z'$ only permutes the values of $\theta'$, also the estimation
	$$
		\int_{\prod_{j\in I} [0,1]^{M_j}}
		(
			f((x_{H,B\setminus\{i\}}(t',\theta'+z')))
			-g_{H,I,t'}(\theta'+z')
		)^2 
		\deld\lebesque^{\coprod_{j\in I} M_j}(\theta')<\bar \gamma^2
	$$
	for all $z'$. However, the function $g_{H,I,t'}(\theta'+z')$ by design only depends on some coordinates and especially not on $z_i$ as $i\notin I$.
	Thus, 
	\begin{align*}
		D_i(t,\theta)
		&\coloneqq 
		\sqrt{\int_{\IT^B} (H_w(t,\theta,z) - \IE_{z_i}H_w(t,\theta,z))^2\deld \lebesque^{B}(z)}
		\\&\le 
		\sqrt{\int_{\IT^B} (H_w(t,\theta,z) -  g_{H,I,t'}(\theta'+z') )^2\deld \lebesque^{B}(z)}
	\end{align*}
	as the arithmetic mean minimizes the quadratic distance.
	Bringing the estimations together, we get for all appropriate $t$ with positive measure that
	\begin{align*}
		&\phantom{{}={}}
		\int_{\prod_{j\in I} [0,1]^{M_j}} 
			D_i(t,\theta)^2
		\deld\lebesque^{\coprod_{j\in I} M_j}(\theta')
		\\&\le
		\int_{\prod_{j\in I} [0,1]^{M_j}} 
			\int_{\IT^B} 
				(H_w(t,\theta,z) -  g_{H,I,t'}(\theta'+z') )^2
			\deld \lebesque^{B}(z)
		\deld\lebesque^{\coprod_{j\in I} M_j}(\theta')
	\\&=
		\int_{\IT^B} 
		\int_{\prod_{j\in I} [0,1]^{M_j}} 
			(H_w(t,\theta,z) -  g_{H,I,t'}(\theta'+z') )^2
		\deld\lebesque^{\coprod_{j\in I} M_j}(\theta')
		\deld \lebesque^{B}(z)	
	\\&\le 
		\int_{\IT^B} 
		2 \cdot \int_{\prod_{j\in I} [0,1]^{M_j}} 
			(H_w(t,\theta,z) -  f((x_{H,B\setminus\{i\}}(t',\theta'+z'))) )^2
		\deld\lebesque^{\coprod_{j\in I} M_j}(\theta')
		\\&\phantom{{}={}}+
		2 \cdot \int_{\prod_{j\in I} [0,1]^{M_j}} 
			(f((x_{H,B\setminus\{i\}}(t',\theta'+z'))) -  g_{H,I,t'}(\theta'+z') )^2
		\deld\lebesque^{\coprod_{j\in I} M_j}(\theta')
		\\&\phantom{{}={}}\deld \lebesque^{B}(z)	
	\\&<
		\int_{\IT^B} 
		2 \cdot \int_{\prod_{j\in I} [0,1]^{M_j}} 
		\bar \gamma^2
		\deld\lebesque^{\coprod_{j\in I} M_j}(\theta')
		+
		2 \cdot \bar \gamma^2
		\deld \lebesque^{B}(z)
	\\&=
		4\bar \gamma^2
	\end{align*}
	using the estimation $2(a-b)^2+2(b-c)^2\ge (a-c)^2$ which is equivalent to $(2b-a-c)^2\ge 0$ and well-known to define $L^2$-norms in the second inequality.
	
	Finally, we can estimate the original set as
	\begin{align*}
			&\cdual\nu_w \{(x_w(t,\theta),t,\theta) 
			\mid \exists i,j\in B\colon t_j=1 \land t_i=0 \land D_i(t,\theta)\ge \beta \}
		\\&\le 
			\sum_{i\in B}\cdual\nu_w \{(x_w(t,\theta),t,\theta) 
			\mid \exists j\in B\colon t_j=1 \land t_i=0 \land D_i(t,\theta)^2/ \beta^2\ge 1 \}
		\\&\le 
			\sum_{i\in B}
			\int_{\{(x_w(t,\theta),t,\theta) 
				\mid \exists j\in B\colon t_j=1 \land t_i=0 \}}
				D_i(t,\theta)^2/ \beta^2
			\deld \cdual\nu_w (x_w(t,\theta),t,\theta)
		\\&<
		\frac{4\bar \gamma^2 \bar s}{\bar \beta ^2} 
		\le \bar \varepsilon
	\end{align*}
	showing the claim.
\end{proof}

\begin{lemma} \label{Lemma_ApproximateWithLambda}
	Consider $w\in W_+$. Now,
	$$ 
		\cdual\nu_w \{(x_w(t,\theta),t,\theta) 
		\mid \exists i\in B\colon t_i\ne 0 \land D_i(t,\theta)\ge \beta \}
		\le \bar \lambda \bar K<\bar \varepsilon
	$$
	holds.
\end{lemma}
\begin{proof}
	Essentially $\bar \lambda$ was chosen large enough to ensure this probability and additionally using the $\bar K$-bound from Lemma~\ref{Lemma_SecondJuntaApplied}.
	
	We fix any $i$ and $t$ with $t_i\ne 0$ and define $t'\in [0,1]^B$ to be the vector $(0,\dots,0,-1/\bar m,0,\dots,0)$ with the nonzero entry at position $i$. Recall that the measure $\rho$ relates $t$ and $t+t'$ with a factor of $\lambda$. As $H$ is $\bar L$-Lipschitz continuous in $t$, it follows
	\begin{align*}
		&|D_i(t,\theta)-D_i(t+t',\theta)|
		\\&= %start of the two-lines
		\left|\sqrt{\int_{\IT^B} (H_w(t,\theta,z) - \IE_{z_i}H_w(t,\theta,z))^2\deld \lebesque^{B}(z)} \right.
		\\&\phantom{{}={}}
		\left.- \sqrt{\int_{\IT^B} (H_w(t+t',\theta,z) - \IE_{z_i}H_w(t+t',\theta,z))^2\deld \lebesque^{B}(z)}\,\right|
		\\&\le %start of the two-lines
		\left|\sqrt{\int_{\IT^B} (H_w(t,\theta,z) - H_w(t+t',\theta,z) )^2\deld \lebesque^{B}(z)} \right.
		\\&\phantom{{}={}}
		\left.+ \sqrt{\int_{\IT^B} ( \IE_{z_i}H_w(t,\theta,z)- \IE_{z_i}H_w(t+t',\theta,z))^2\deld \lebesque^{B}(z)}\,\right|
		\\&\le %start of line
		\left|\sqrt{\int_{\IT^B} (L|t'| )^2\deld \lebesque^{B}(z)} + \sqrt{\int_{\IT^B} (L|t'|)^2\deld \lebesque^{B}(z)}\,\right|
		\\&=
		2L|t'|= 2L/m<\beta/2
	\end{align*}
	and 
	\begin{align*}
		D_i(t,\theta)\ge \beta \implies D_i(t+t',\theta)\ge \beta/2.
	\end{align*}
	Summing up over all $i$, we get
	\begin{align*}
		&\cdual\nu_w \{(x_w(t,\theta),t,\theta) 
		\mid \exists i\in B\colon t_i\ne 0 \land D_i(t,\theta)\ge \beta \}
		\\&\le 
		\sum_{i\in B} \cdual\nu_w \{(x_w(t,\theta),t,\theta) 
		\mid t_i\ne 0 \land D_i(t,\theta)\ge \beta \}
		\\&\le 
		\sum_{i\in B} \cdual\nu_w \{(x_w(t,\theta),t,\theta) 
		\mid t_i\ne 0 \land D_i(t+t',\theta)\ge \beta/2 \}
		\\&=
		\sum_{i\in B} \bar \lambda \cdot \cdual\nu_w \{(x_w(t,\theta),t,\theta) 
		\mid t_i\ne 1 \land D_i(t,\theta)\ge \beta/2 \}
		\\&\le 
		\bar \lambda \cdot \sum_{i\in B} \cdual\nu_w \{(x_w(t,\theta),t,\theta) 
		\mid D_i(t,\theta)\ge \beta/2 \}
		\\&\le \bar \lambda \bar K<\bar \varepsilon
	\end{align*}
	where the last inequality uses Lemma~\ref{Lemma_SecondJuntaApplied}.
\end{proof}
\begin{corollary} \label{Corollary_BoundZ}
	Consider $w\in W_+$. Now,
	$$ 
	\cdual\nu_w \{(x_w(t,\theta),t,\theta) 
	\mid \exists j\in B\colon t_j=1 \land
	\IE_z H_w(t,\theta,z)^2>\bar\varepsilon
	\}< 2 \bar \varepsilon
	$$
	holds.
\end{corollary}
\begin{proof}
	If $\IE_z H_w(t,\theta,z)^2>\bar\varepsilon$, then there is $i$ with $D_i(t,\theta)\ge \bar \beta$ by Lemma~\ref{Lemma_SecondJuntaApplied}. However, this case was bounded by Lemma~\ref{Lemma_UsingAllinged} and Lemma~\ref{Lemma_ApproximateWithLambda} by $\bar \varepsilon$ each.
\end{proof}
\begin{claim} \label{Claim_greatHw}
	Consider $w\in W_+$. Now,
	$$
		\int_{\{(t,\theta)\in [0,1]^I \times [0,1]^{\prod_{l\in I} M_l} \mid \exists j\in B\colon t_j=1\}} \IE_z (H_w(t,\theta,z))^2
		\deld (\rho^B \times \lebesque^{\prod_{j\in B} M_j}) (t,\theta)
	$$
	is bounded by $<3\bar\varepsilon$ from above.
\end{claim}
\begin{proof}
	Since the function is bounded by $1$, which is also a bound for the measure of the total space, its integral is surely smaller than 
	$$
		\bar \varepsilon \cdot 1 + 1\cdot (\rho^B \times \lebesque^{\prod_{j\in B} M_j})(\{(t,\theta )\mid \IE_z (H_w(t,\theta,z))^2>\bar\varepsilon\})< 3\bar\varepsilon
	$$
	by Corollary~\ref{Corollary_BoundZ}.
\end{proof}
\begin{claim}[Independence of $z$] \label{Claim_zIndependence}
	Consider $w\in W_+$ and $t\in [0,1]^B$ with $j\in B, t_j=1$.
	Then, 
	\begin{align*}
		&\int_{[0,1]^{\prod_{l\in I} M_l} } \IE_z (H_w(t,\theta,z))^2
		\deld \lebesque^{\prod_{l\in B} M_l} (\theta)
		\\= {}&\int_{[0,1]^{\prod_{l\in I} M_l} } (H_w(t,\theta,0))^2
		\deld \lebesque^{\prod_{l\in B} M_l} (\theta)
	\end{align*}
	holds.
\end{claim}
\begin{proof}
	Hide $z$ in $\theta_{j} \in \IT^{M_j}$ by a substitution.
\end{proof}
\begin{lemma} \label{lemma_continuousBound}
	Consider $w\in W_+$. Now,
	$$
		\int_{M_w \times [0,1]^B \times [0,1]^{\prod_{j\in B} M_j}} 
		(f(x))^2
		\deld \cdual \nu_w(x,t,\theta)< 4\bar \varepsilon
	$$
	holds.
\end{lemma}
\begin{proof}
	Excluding $t$ without any $t_j=1$ excludes a set with measure $<  \bar\varepsilon$ by Lemma~\ref{lemma_fullCoefficientNeeded}. As $f(x)^2\le 1$, this summand of the integral is bounded by $\bar\varepsilon$.
	
	For other $t$, we have $(f(x))^2=H_w(t,\theta,0)$ whenever $\nu_w$ is positive and can rewrite the function using Claim~\ref{Claim_zIndependence} into a function using $z$ that was already bounded by $3\bar \varepsilon$ in Claim~\ref{Claim_greatHw}.
\end{proof}
We can discretize this result and omit the assumption of an $j\in B$ with $t_j=1$ by adding two more $\bar\varepsilon$.
\begin{lemma}\label{lemma_fToEpsilon}
	Consider $w\in W_+$. Now,
	$$
	\int_{V(\ix)} (f(x))^2
	\deld \nu_w(x)<5\bar \varepsilon
	$$
	holds.
\end{lemma}
\begin{proof}
	By definition, the integral equals
	$$
		\int_{M_w \times [0,1]^B \times [0,1]^{\prod_{j\in B} M_j}} (f(x))^2
		\deld \fdual \nu_w(x,t,\theta).
	$$
	 Now, we can apply Lemma~\ref{lemma_continuousBound} and have to show that the difference of the two integrals with different measures is bounded by $\bar \varepsilon$. 
	 However, this is true as the two measures $\cdual\nu_w$ and $\fdual\nu_w$ differ in the way that $\cdual\nu_w$ considers $\theta$ continuously while $\fdual\nu_w$ considers $\theta$ as discrete samples. Define $\fdual\theta$ component wise as 
	 $$
	 	\fdual\theta_{j,l}=\frac{\lfloor \theta_{j,l} P \rfloor}{P}.
	 $$
	 Then, we get that
	 
	 \begin{align*}
	 	&\phantom{{}={}}
	 	\int_X (f(x))^2
	 	\deld \fdual \nu_w(x,t,\theta)
	 	- \int_X (f(x))^2
	 	\deld \cdual \nu_w(x,t,\theta)
	 	\\
	 	&= 
	 	\int_X
	 	(f(x_w(t,\theta)))^2 - (f(x_w(t,\fdual\theta)))^2
	 	\deld \cdual \nu_w(x,t,\theta)
	 	\\
	 	&=\int_X 
	 	(f(x_w(t,\theta))+f(x_w(t,\fdual\theta)))
	 	(f(x_w(t,\theta))-f(x_w(t,\fdual\theta)))
	 	\deld \cdual \nu_w(x,t,\theta)
	 	\\
	 	&\le \int_X 2 \cdot \bar L\bar s|\theta-\fdual\theta| \deld \cdual \nu_w(x,t,\theta)
	 	\\
	 	&\le 2 \bar L\bar s / \bar P \int_X 1 \deld \cdual \nu_w(x,t,\theta)
	 	\\
	 	&= 2 \bar L\bar s / \bar P <\bar\varepsilon
	 \end{align*}
	 where all integrals are over the area $X=M_w \times [0,1]^B \times [0,1]^{\prod_{j\in B} M_j}$. The first inequality uses $f(x)\le 1$, the $\bar L$-Lipschitz continuity of $f$ and the $\bar s$-Lipschitz continuity of $x_w$. The second inequality bounds the step size by $1/\bar P$.
\end{proof}

\begin{lemma}[c.f. Proposition 6.1 in \cite{OAI_PCSP}] \label{Lemma_AnalysisResult}
	Consider $w\in W_+$. Now,
	$
		\nu_w(\mathcal{A})<5/2\cdot  \bar \varepsilon.
	$
\end{lemma}
\begin{proof}
	$f(x)^2=1$ whenever $x\in \mathcal{A}$ or $T(x)\in \mathcal{A}$. These two cases exclude each other and as $f(x)\le 1$, it follows
	$$2 \nu_w(\mathcal{A})= \nu_w(\mathcal{A} \cup T(\mathcal{A})) \le \int_{V(\ix)} (f(x))^2
	\deld \nu_w(x)<5\bar \varepsilon$$ by Lemma~\ref{lemma_fToEpsilon}.
\end{proof}
\subsection{Combining Estimations}
\label{Section_corr_Combined}
We can now combine the estimations.
\begin{proposition} \label{Prop_MaxMeasure}
	Assume that $\inst$ is a negative instance of the $\bar \sigma$-gap $\bar A$-multi-label cover problem with bound $B$ and that $\mathcal{A}$ is an independent vertex set in $\ix$. Then, it follows that $\nu(\mathcal{A})<5\bar\varepsilon<1/k$.
\end{proposition}
\begin{proof}
	We have
	\begin{align*}
		&\nu(\mathcal A)
		\\&= \sum_{w\in W_- \text{ not aligned}} \mu_W(w) \nu_w\{\mathcal A\}
		+\sum_{w\in W_- \text{ aligned}} \mu_W(w) \nu_w\{\mathcal A\}
		\\&\le  \sum_{w\in W_- \text{ not aligned}}  \mu_W(w) 
		+ \max \{\nu_w\{\mathcal A\} \mid w\in W_- \text{aligned}  \}
		\\&< 2\bar \varepsilon + 5/2 \bar \varepsilon <5 \bar \varepsilon
	\end{align*}
	using Lemma~\ref{Lemma_notAlligned} and Lemma~\ref{Lemma_AnalysisResult}.
\end{proof}

\begin{proofof}{Theorem~\ref{TheoremReduction1}}
	It is easy to see that the construction is in polynomial time (using that $k$ and $B$ are fixed). The first claim is shown in Lemma~\ref{lemma_easyDirection}. For the second claim, we use the 
	notation from Definition~\ref{algdef_constants} and~\ref{algdef_measures}.
	Chose $\sigma=\bar \sigma, A=\bar A$. 
	Now, $G$ and $\ix=\ix_{k,\bar D,\bar r}$ are homomorphically equivalent (Lemma~\ref{lemma_homequiv}, thus it suffices to show that there is no homomorphism $\ix\to K_k$. 
	Assume to the contrary that there would be a homomorphism. Every preimage of the $k$ points in $K_k$ is an independent set and thus, by Proposition~\ref{Prop_MaxMeasure}, has $\nu$-measure at most $5\bar\varepsilon$.
	Thus, $\xi$ as the full preimage of $K_k$ would have measure at most $k\cdot 5\bar \varepsilon<1$, a contradiction.
\end{proofof}

\begin{proofof}{Theorem~\ref{TheoremMain}}
	If $G$ is bipartite then already $G$-coloring can be done in polynomial time without a promise: If $G$ has no vertices or no edges this problem is trivial. Otherwise $G$-coloring is exactly 2-coloring which can be done in polynomial time as well. If $H$ has a loop, then the $H$-coloring problem is trivial as every vertex of the input graph can be mapped to this loop. If $H$ is bipartite or $G$ has a loop then the same property follows for the other graph as well. Thus, in those cases, the promise graph coloring problem can be solved in polynomial time.
	
	For general $G$ and $H$, the promise graph coloring is in NP as it can be checked in polynomial time that a given $H$-coloring works if there exists one.
	The promise graph coloring problem reduces in polynomial time to the gap multi-label cover problem with a bound by Theorem~\ref{TheoremReduction1} which is NP-hard by Corollary~\ref{Corollary_BoundLabelCover} if $G$ and $H$ are non-bipartite and without loop. Thus, the problem is NP complete in this case.
\end{proofof}

\appendix
\section{Generalizations}

Theorem~\ref{TheoremMain} completely solves promise graph coloring for undirected graphs. However, we may define similar decision problems and may wonder for which problems we are able to show hardness in a similar way.

\begin{theorem}
	For every (finite) non-bipartite graph $G$ without loop and every integer $k\in (0,1)$, the following promise decision problem is NP-hard.
	\begin{description}
		\item[Input:] A graph $K$,
		\item[Output: Yes] if there is a graph homomorphism $K\to G$,
		\item[Output: No] if there is no independent set in $K$ with size $V(G)/k$.
	\end{description}
	In other cases, there can be any output.
	
\end{theorem}
\begin{proof}
	We use a similar argument as in the proofs of Theorem~\ref{TheoremReduction1} and Theorem~\ref{TheoremMain} with the graph $\ix'$ from Remark~\ref{remark_cliqueResult} instead of the graph $G$. The graph $\ix'$ can be constructed in polynomial time from a gap multi-label cover with bounded instance. It defines a reduction from gap multi-label cover with bound to the above problem as there is a graph homomorphism $\ix'\to \ix$ thus positive instances induce positive instances and Proposition~\ref{Prop_MaxMeasure} shows that there is no independent set in $\ix$ (with respect to the counting measure on $\ix$) containing a fraction of $5\bar\varepsilon<1/k$ of the vertices.
\end{proof}

Instead of undirected graphs, we may also consider directed graphs or relational structures. As the main ingredient to the proof was the Junta theorem idea, we may guess that a similar reduction works whenever we have a similar theorem.
\begin{conjecture}
	Let $V$ be a (finite) set, $n$ a positive integer and $E\subset V^n$ a subset.
	
	A subset $I\subseteq V^m$ is called independent, if there is no matrix with columns in $I$ and rows in $E$. 
	
	Assume that there is a probability measure $\mu$ on $V$ and a function $j\colon (0,1)\to \IN$ such that for every $\varepsilon\in(0,1)$ and every independent set $I\subset V^m$ there is a set $J\subset V^m$, depending on at most $j$ coordinates such that $\mu(I\setminus J)\le \varepsilon$ and such that the number of different matrices with columns in $J$ and rows in 
	$E$ is at most $\varepsilon\cdot |E|$.
	
	In this case, the promise problem
	\begin{description}
		\item[Input:] A set $K$ with a subset $L$ of $K^m$,
		\item[Output: Yes] if there is a map $K\to V$ that maps tuples in $L$ to tuples in $E$,
		\item[Output: No] if there is no independent set in $K$ with size $V(G)/k$.
	\end{description}
	is NP complete.
\end{conjecture}

\bibliography{Quellen} 
\bibliographystyle{alpha}

\end{document}